\documentclass[10pt]{amsart}

\usepackage{latexsym,amssymb,amsmath,graphics}
\usepackage[dvips]{graphicx}

\def\a{\alpha}
\def\aa{{\mathcal A}}
\def\bb{{\mathcal B}}
\def\cc{{\mathcal C}}
\def\dd{{\mathcal D}}
\def\ee{{\mathcal E}}
\def\ff{{\mathcal F}}
\def\gg{{\mathcal G}}
\def\hh{{\mathcal H}}
\def\jj{{\mathcal J}}
\def\kk{{\mathcal K}}
\def\lll{{\mathcal L}}
\def\mm{{\mathcal M}}
\def\nn{{\mathcal N}}
\def\pp{{\mathcal P}}
\def\qq{{\mathcal Q}}
\def\rr{{\mathcal R}}
\def\ss{{\mathcal S}}
\def\ttt{{\mathcal T}}
\def\uu{{\mathcal U}}
\def\ww{{\mathcal W}}
\def\xx{{\mathcal X}}
\def\yy{{\mathcal Y}}
\def\zzz{{\mathcal Z}}

\def\8{\infty }
\def\ie{i.e.\ }
\def\tq{\ :\ }
\def\ffi{\varphi}
\def\eps{\varepsilon}
\def\dst{\displaystyle}
\def\scp{\scriptstyle}
\def\eg{{\it e.g.\ }}
\def\etc{{\it etc}}
\def\supp{{\mathrm{supp}\,}}
\def\ovl{\overline}

\def\A{{\mathbb{A}}}
\def\B{{\mathbb{B}}}
\def\C{{\mathbb{C}}}
\def\hil{{\mathbb{H}}}
\def\M{{\mathbb{M}}}
\def\N{{\mathbb{N}}}
\def\P{{\mathbb{P}}}
\def\Q{{\mathbb{Q}}}
\def\R{{\mathbb{R}}}
\def\S{{\mathbb{S}}}
\def\T{{\mathbb{T}}}
\def\X{{\mathbb{X}}}
\def\Y{{\mathbb{Y}}}
\def\Z{{\mathbb{Z}}}

\def\w{{\mathrm{w}}}

\def\hei{{\mathbf{H}^n}}

\newcommand{\norm}[1]{{\left\|{#1}\right\|}}
\newcommand{\ent}[1]{{\left[{#1}\right]}}
\newcommand{\abs}[1]{{\left|{#1}\right|}}
\newcommand{\scal}[1]{{\left\langle{#1}\right\rangle}}

\newenvironment{notation}[1][]{\vskip1pt\noindent\rm\textit{Notation.}\
}{\rm\vskip1pt}
\newenvironment{remark}[1][]{\vskip1pt\noindent\rm\textit{Remark.}\
}{\rm\vskip1pt}
\newenvironment{definition}[1][]{\vskip3pt\noindent\sl\textbf{Definition.}\
}{\rm\vskip3pt}

\renewcommand{\Im}{\mathrm{Im}\,}
\renewcommand{\Re}{\mathrm{Re}\,}

\newtheorem{lemma}{Lemma}[section]
\newtheorem{proposition}[lemma]{Proposition}
\newtheorem{theorem}[lemma]{Theorem}
\newtheorem{mth}[lemma]{Main Theorem}

\renewcommand{\theequation}{\thesection.\arabic{equation}}
\begin{document}

\title[Bounded Hua-harmonic functions]{Maximum boundary regularity of
bounded Hua-harmonic functions on tube domains}
\author[Aline Bonami, Dariusz Buraczewski, Ewa Damek, Andrzej
Hulanicki \& Philippe Jaming]{Aline Bonami, Dariusz Buraczewski$^1$, Ewa Damek$^1$, Andrzej
Hulanicki$^1$ \& Philippe Jaming}
\address{A. Bonami \& P. Jaming\\ Universit\'e d'Orl\'eans\\
Facult\'e des Sciences\\
D\'epartement de Ma\-th\'e\-ma\-ti\-ques\\BP 6759\\ F 45067 Orl\'eans
Cedex 2\\
France}
\email{bonami@labomath.univ-orleans.fr\\
jaming@labomath.univ-orleans.fr}

\address{D. Buraczewski, E. Damek \& A. Hulanicki\\ 50-384 Wroclaw\\
pl. Grunwaldzki 2/4\\ Poland}
\email{dbura@math.uni.wroc.pl\\ edamek@math.uni.wroc.pl\\
hulanick@math.uni.wroc.pl}

\begin{abstract}
In this paper we prove that bounded Hua-harmonic functions on tube domains that satisfy some boundary regularity condition are necessarily pluriharmonic. In doing so, we show that a similar theorem is true on one-dimensional extensions of the Heisenberg group or equivalently on the Siegel upper half-plane.
\end{abstract}


\thanks{ Research partially financed by the {\it European Commission}
IHP Network 2002-2006 {\it Harmonic Analysis and Related Problems}
(Contract Number: HPRN-CT-2001-00273 - HARP).\\
\indent
$^1$The authors were partly supported by KBN grant 5 PO3A 02821
and Foundation for Polish Sciences Subsidy 3/99}

\maketitle
\section{Introduction}

Let $\uu^n\subset\C^{n+1}$ be the  Siegel upper half-plane,
defined by
$$
\uu^n=\left\{z\in\C^{n+1}\,:\ \Im z_{n+1}>\sum_{j=1}^n\abs{z_j}^2\right\}.
$$
Let $F$ be the Poisson-Szeg\"{o} integral of some boundary function
$f$. It has been known for a long time now that $F$, which is
smooth inside $\uu^n$, cannot have bounded transversal Euclidean
derivatives up to the boundary, unless $F$ is a pluriharmonic
function. More precisely, in $\uu^n$ we introduce coordinates
$\zeta=(z_1,...,z_n)$, $z_{n+1}= t+i|\zeta|^2+ia$, with $t\in \R$,
and $a>0$, so that, for the complex structure $J$, we have
$J\partial_a=-\partial_t$. Then, if $F$ is harmonic with respect
to the invariant Laplacian, the derivatives $\partial_a^kF(\zeta,
t+i|\zeta|^2+ia)$ cannot be bounded for large $k$ unless $F$ is
pluriharmonic. Such results may be found in the work of Graham
(see \cite{Gr1} and \cite{Gr2}). The conditions are easier to
describe for functions on the unit ball in $\C^{n+1}$, that is,
for the bounded realization of $\uu^n$, see \cite{BBG}. In both
cases, a central role is played by the invariant Laplacian $L$.
Indeed, Poisson-Szeg\"{o} integrals are characterized by the fact that
they are annihilated by $L$ so that these results are merely
results about $L$-harmonic functions.

\smallskip

The aim of this paper is to show that this property holds in a general
context. More precisely, we consider an irreducible symmetric
Siegel domain of tube type which may be written as
$$
\dd=V+i\Omega\;\subset V^{\C},
$$
where $V$ is a real Euclidean space of dimension $m$ and $\Omega$
is an irreducible symmetric cone inside $V$. Typical examples  are
given when one chooses for $\Omega$ the forward light cone or the
cone of positive definite matrices. For such domains, Poisson-Szeg\"{o}
integrals have been characterized in terms of
differential operators of order 2. Let us give more details, and
introduce the Hua system which plays the same role as the
invariant Laplacian in the case of $\uu^n$.
This system can be defined geometrically for a more general
domain $\dd $ in $\C^m$ which is
biholomorphically equivalent to a bounded domain.  For more
details we refer to  \cite{BBDHPT} where the following
construction, which is inspired by the work of Wallach, is
thoroughly discussed.

Let $T^{\C}$ be the
complexified tangent bundle of the complex domain $\dd$, let
$\jj$ be the complex structure, and let $T^{1,0}$
and $T^{0,1}$ be the eigenspaces of $\jj $ such that $\jj
|_{T^{1,0}}= i \text{Id}$, $\jj |_{T^{0,1}}=-i\text{Id}$.
The Riemannian connection $\nabla
$ induced by the Bergman metric on $\dd$ preserves $T^{1,0}(\dd )$
and so does the curvature tensor $R(Z,W)$. For $Z,W$ two complex
vector fields we denote by
  $R(Z,W)=\nabla _Z\nabla _W-\nabla _W\nabla _Z-\nabla _{[Z,W]}$
the  curvature tensor restricted to $T^{1,0}(\dd )$. For $f$  a
smooth function on $\dd $, let
\begin{equation}
\Delta (Z,W)f=(Z\overline{W}-\nabla_{Z}\overline{W})f=
(\overline{W}Z-\nabla_{\overline{W}}Z)f\,.
\label{delta}
\end{equation}
Then $\Delta (Z,W)$ annihilates both holomorphic and
anti-holomorphic functions, and consequently, the pluriharmonic
functions.

Given an orthonormal frame $E_1,E_2,\cdots,E_m$ of $T^{1,0}(\dd)$
for the canonical Hermitian product associated to the Bergman metric,
the Hua system is defined as follows:
\begin{equation}
\hil f=\sum _{j,k}\bigl(\Delta (E_j, E_k)f\bigr) R(\overline
E_{j}, E_{k}) \label{Hua-inv}\,.
\end{equation}
The Hua system does not depend on the choice of the orthonormal
frame and it is invariant with respect to biholomorphisms. By
definition, $\hil$-harmonic functions are functions  which are
annihilated by $\hil$.

When $\dd$ is a symmetric Siegel domain of
tube type, this system is known to characterize the Poisson-Szeg\"o
integrals (see \cite {FK} and \cite {JK}). This means that a
function on $\dd$ is $\hil$-harmonic if, and only if, it is the
Poisson-Szeg\"o integral of a hyperfunction on the Shilov
boundary.

Our main theorem, may then be stated as follows:

\smallskip

{\bf Main Theorem.} {\it Let $\dd$ be an irreducible symmetric
domain of tube type. There exists $k$ (depending on the dimension
and the rank) such that, if $F$ is a bounded $\hil$-harmonic and
has bounded derivatives up to the order $k$, then $F$ is
pluriharmonic.}

\smallskip

More precise and weaker conditions are given in the sequel of the
paper. They are linked to the description of the domain in terms
of a solvable group of linear automorphisms and may be expressed
in the distribution sense. Roughly speaking, we show that there
exists an open dense subset $\widetilde{\partial\dd}$ of $\partial
\dd$, such that for every point $p\in \widetilde{\partial\dd}$
there is an open  neighborhood of $p$ in $V+iV$, called $U$, and a
smooth foliation  $U_a,\ a\in [0,\epsilon ]$, of $U\cap
\overline{\dd}$ such that $U_0=\partial\dd\cap U$ and
$J\partial_a$ is tangential to $U_a$. Moreover,
$\widetilde{\partial\dd}$ is invariant under a group of
biholomorphisms of $\dd$ acting transitively on
 the family of foliations. Locally, the parameter $a$ plays the same role as the coordinate $a$ in
  the Siegel upper half-plane. If
$F$ is Hua-harmonic and for large $k$, $\partial_a^kF(\cdot, a)$
is bounded as $a\to 0$ for all such parameters $a$, then we prove
that $F$ is pluriharmonic.

\smallskip

Let us remark that many sufficient conditions
for pluriharmonicity, which can be written
 in terms of families of second order operators, have been given by some of the authors
(see \cite{BDH} and \cite{DHMP}). We rely deeply on this previous
work. Another main tool of the proof is the fact that we can
reduce to the same kind of problem on the domain $\uu^n$ for a
variant of the invariant Laplacian $L$. So, we are lead to
consider a whole family of second order operators on $\uu^n$, and we
give sufficient conditions so that functions which are annihilated
by them are pluriharmonic. A key tool for this is the existence of
a boundary equation which implies pluriharmonicity.

\smallskip

One may ask whether this type of results can be obtained in the
more general setting of irreducible symmetric Siegel domains,
which contains both those  of tube type as well as the Siegel
upper half-plane $\uu^n$, which is of type II. Unfortunately, much
less is known for higher rank and type II. In particular, for the
Hua system that we described, $\hil$-harmonic functions are
already pluriharmonic (see \cite{BBDHPT} and \cite{B}), and such tools
are missing for studying Poisson-Szeg\"{o} integrals (see
\cite{BBDHPT} for more comments). So, the kind of proof that we
are giving here cannot be adapted to such a general situation.

\smallskip

Let us mention that the same phenomenon has also been studied by
B. Trojan, using a direct computation of the Poisson-Szeg\"{o} kernel
and looking at its singularities (see \cite{Tr}).

\smallskip

The article is organized as follows. In the next section, we
introduce all necessary notations on the Heisenberg group and its
one-dimensional extensions. We then prove a characterization of
bounded pluriharmonic functions in terms of boundary conditions
and give sufficient regularity conditions for our harmonic
functions which imply pluriharmonicity. In the next section,
 we give all preliminaries that we need on Siegel
domains of tube type, including a precise expression of Hua
operators in terms of the description of the cone using Jordan
algebras. Finally we prove that Hua-harmonic operators have a
maximum boundary regularity, unless they are pluriharmonic.

Some technical results on special functions, as well as another
proof of the main theorem, are given in the appendix.

\section{The theorem on the Heisenberg group}
\label{heisenberg}\setcounter{equation}{0}

\subsection{Preliminaries on the Heisenberg group}

In this section we  recall some notations and results on the
Heisenberg group. Our basic
reference is \cite{St}. We will try to keep self-contained up
to some results in the two last chapters of that book.

The Heisenberg group $\hei$ is the set
$$
\hei=\C^n\times\R=\{[\zeta,t]\,:\ \zeta\in\C,\ t\in\R\}
$$
endowed with the multiplication law
$$
[\zeta,t].[\eta,s]=[\zeta+\eta,t+s+2\Im(\zeta\overline{\eta})].
$$
We  write $[\zeta,t]$ for a typical element of $\hei$ and
$\zeta=(\zeta_1,\ldots,\zeta_n)=(x_1+iy_1,\ldots,x_n+iy_n)\in\C^n$.

The Euclidean Lebesgue measure $dx\, dy\,dt$, which we will note
as well $d\zeta\,dt$ or $d\omega$ if $\omega= [\zeta,t]$, is both
left and right invariant under the action of $\hei$. Convolution
on $\hei$ is then given by
$$
f*g([\zeta,t])=\int_{\hei}f([\eta,s])g([\eta,s]^{-1}[\zeta,t])d\eta
ds.
$$

We next consider the domain $\uu^n\subset\C^{n+1}$ and its boundary
$b\uu^n$ defined by
$$
\uu^n=\left\{z\in\C^{n+1}\,:\ \Im
z_{n+1}>\sum_{j=1}^n\abs{z_j}^2\right\},\quad
b\uu^n=\left\{z\in\C^{n+1}\,:\ \Im
z_{n+1}=\sum_{j=1}^n\abs{z_j}^2\right\}.
$$
A typical element of $\uu^n$ or $b\uu^n$ is denoted by
$z=(z',z_{n+1})$.
The Heisenberg group $\hei$ acts on $\uu^n$ and $b\uu^n$ by
$$
[\zeta,t](z',z_{n+1})=(z'+\zeta,z_{n+1}+t+2iz'\overline{\zeta}+i\abs{\zeta}^2).
$$
Moreover, this action is simply transitive on the boundary $b\uu^n$
which allows us to identify
elements of $\hei$ with elements of $b\uu^n$ by the action of $\hei$
on the origin $(0,0)$~:
$$
\hei\ni[\zeta,t]\mapsto [\zeta,t](0,0)=(\zeta,t+i\abs{\zeta}^2)\in
b\uu^n.
$$

Further, we write $r(z)=\Im(z_{n+1})-\abs{z'}^2$, so that
$\uu^n=\{z\,:\ r(z)>0\}$ and $b\uu^n=\{z\,:\ r(z)=0\}$.
The Heisenberg group $\hei$ acts simply transitively on each level set so that
each one of them can be
identified with $\hei$.

Next, we  consider on $\hei$ the  left-invariant vector
fields
$$
X_j=\frac{\partial}{\partial x_j}+2y_j\frac{\partial}{\partial
t},\quad
Y_j=\frac{\partial}{\partial y_j}-2x_j\frac{\partial}{\partial
t},\quad
T=\frac{\partial}{\partial t}.
$$
One has $[Y_j,X_k]=4\delta_{j,k}T$ while all other
commutators vanish.

We  define the holomorphic and anti-holomorphic left invariant vector fields
$$
\overline{Z}_j=\frac{1}{2}(X_j+iY_j)=\frac{\partial}{\partial\overline{\zeta}_j}-i\zeta_j\frac{\partial}{\partial
t},\quad
Z_j=\frac{1}{2}(X_j-iY_j)=\frac{\partial}{\partial\overline{\zeta}_j}+i\overline{\zeta}_j\frac{\partial}{\partial
t}
$$
so that $[\overline{Z}_j,Z_k]=2i\delta_{j,k}T$ while all other
commutators are $0$.

The following operators then play a fundamental role in complex
analysis on $\uu^n$
(see \cite{St}):
$$
\lll_\alpha=-\frac{1}{2}\sum_{j=1}^n(Z_j\overline{Z}_j+\overline{Z}_jZ_j)+i\alpha
T
=-\frac{1}{4}\sum_{j=1}^n(X_j^2+Y_j^2)+i\alpha T,
$$
where $\alpha\in\R$ is a parameter. We will simply write $\lll$ for
$\lll_0$.

The Cauchy kernel on $\hei$ is then given, for $[\zeta,t]\not=[0,0]$,
by $K([\zeta,t])=c(t+i\abs{\zeta}^2)^{-n-1}$
with $c=2^{n-1}i^{n+1}n!/\pi^{n+1}$ and the Cauchy-Szeg\"o projection
is defined by
$C(f)=f*K$, where $K$ defines as well the principal value
distribution which is associated to this kernel.

Define
$$
\Phi([\zeta,t])=\frac{2^{n-2}(n-1)!}{\pi^{n+1}}\log\left(\frac{\abs{\zeta}^2-it}{\abs{\zeta}^2+it}\right)(\abs{\zeta}^2-it)^{-n}
$$
where
$$
\log\left(\frac{\abs{\zeta}^2-it}{\abs{\zeta}^2+it}\right)
=\log(\abs{\zeta}^2-it)-\log(\abs{\zeta}^2+it)
$$
and the logarithms are taken to be their principal branch in the
right half-plane.
Write $\widetilde{S}(f)=f*\Phi$. From Formula (53) page 616 in
\cite{St}, we then know that if $f$ is a smooth compactly supported
function, then
\begin{equation}
\label{Cauchy}
\lll_n\widetilde{S}(f)=\widetilde{S}(\lll_nf)=f-C(f).
\end{equation}

\subsection{One dimensional extension of the Heisenberg group}

We now define the non-isotropic dilations on $\uu^n$ and $b\uu^n$:
for $\delta>0$ and $z=(z',z_{n+1})\in\uu^n$ or $b\uu^n$, we write
$$
\delta\circ z=(\delta z',\delta^2 z_{n+1}).
$$
These dilations preserve $\uu^n$ and $b\uu^n$
whereas
$r(\delta\circ z)=\delta^2r(z)$.

We now consider the semi-direct extension of $\hei$, $S=\hei\R^+_*$
where the action of $\R^+_*$
on $\hei$ is given by $a[\zeta,t]=[a^{1/2}\zeta,at]$. We will write
$[\zeta,t,a]=\bigl[[\zeta,t],a\bigr]$ for a typical element of $S$
with $[\zeta,t]\in\hei$ and $a\in\R^+_*$. The group law of $S$ is
then given by
$$
[\zeta,t,a][\eta,s,b]=\bigl[[\zeta,t][a^{1/2}\eta,as],ab\bigr]=[\zeta+a^{1/2}\eta,t+as+2a^{1/2}
\Im(\overline{\zeta}\eta),ab].
$$
We then extend the definition of the previous vector fields from
$\hei$ to $S$ so that the following are left-invariant vector
fields on $S$:
$$
a^{1/2}X_j,\ a^{1/2}Y_j,\ aT,\ a\partial_a=a\frac{\partial}{\partial
a},\ a^{1/2}Z_j,\
a^{1/2}\overline{Z}_j.
$$
We next define $Z_{n+1}=\frac{1}{2}(T-i\partial_a)$ and
$\overline{Z}_{n+1}=\frac{1}{2}(T+i\partial_a)$
so that $aZ_{n+1}$ and $a\overline{Z}_{n+1}$ are left-invariant
vector fields on $\uu_n$ (this notation differs from that in \cite{St}
by a harmless factor $2^{1/2}i$).

We  get  $[a^{1/2}X_j,a^{1/2}Y_j]=4\delta_{j,k}aT$,
$[a\partial_a,a^{1/2}X_j]=\frac{1}{2}a^{1/2}X_j$,
$[a\partial_a,a^{1/2}Y_j]=\frac{1}{2}a^{1/2}Y_j$
$[a\partial_a,aT]=aT$, while all other commutators are $0$.

The group $S$ acts simply transitively on $\uu^n$ by
$[0,0,a]z=a^{1/2}\circ z$ and
$[\zeta,t,1]z=[\zeta,t]z$ as defined previously. In other words
$$
[\zeta,t,a]z=[\zeta,t,1][0,0,a]z=[\zeta,t](a^{1/2}z',az_{n+1})=(\zeta+a^{1/2}z',
t+az_{n+1}+2ia^{1/2}z'\overline{\zeta}+i\abs{\zeta}^2).
$$
In particular $[\zeta,t,a](0,i)=(\zeta,t+i\abs{\zeta}^2+ia)$. This
allows us to identify
$S$ with $\uu^n$ and a function $f_s$ on $S$ with a function $f_u$
on $\uu^n$ by
$$
f_s([\zeta,t,a])=f_u(\zeta,t+i\abs{\zeta}^2+ia)\quad\mbox{and}\quad
f_u(z',z_{n+1})=f_s([z',\Re(z_{n+1}),\Im(z_{n+1})-\abs{z'}^2]).
$$
It then follows that, for $k=1,\ldots,n$, $\dst\frac{\partial
f_u}{\partial z_k}=Z_kf_s-2i\overline{z}_kZ_{n+1}f_s$.
We will thus say that
\begin{itemize}
\item[---] $f_s$ is holomorphic if $f_u$ is holomorphic, {\it i.e.}
$\frac{\partial f_u}{\partial\overline{z}_k}=0$ for $k=1,\ldots,n+1$.
This is easily seen to be equivalent to $\overline{Z_k}f_s=0$ for
$k=1,\ldots,n+1$.\\
\item[---] $f_s$ is anti-holomorphic if $f_u$ is anti-holomorphic,
{\it i.e.}
$\frac{\partial f_u}{\partial z_k}=0$ for $k=1,\ldots,n+1$ or
equivalently
if $Z_kf_s=0$ for $k=1,\ldots,n+1$.\\
\item[---] $f_s$ is pluriharmonic if $f_u$ is pluriharmonic, {\it
i.e.}
$\frac{\partial}{\partial
z_k}\frac{\partial}{\partial\overline{z}_j}f_u=0$ for
$j,k=1,\ldots,n+1$.
This is then equivalent to $Z_k\overline{Z}_jf_s=0$ for $1\leq
j\not=k\leq n+1$,
$Z_{n+1}\overline{Z}_{n+1}f_s=0$ and
$(Z_k\overline{Z}_k+2i\overline{Z}_{n+1})f_s=0$ for $k=1,\ldots,n$.
\end{itemize}
In the sequel, we will drop the subscripts $s$ and $u$ for simplicity
as well as the superfluous brackets $[,]$.

\begin{notation} Let us fix $\alpha>0$ and consider the
left-invariant operator $L_\alpha$ given by
$$
L_\alpha=-\alpha a(\lll+n\partial_a)+a^2(\partial_a^2+T^2).
$$
\end{notation}
Of particular interest in the next section will be the case
$\alpha=\frac{1}{2}$.

Let $P_a^\alpha$ be the Poisson kernel for $L_\alpha$, i.e. the
function on $\hei$ that establishes a one-to-one correspondence
between bounded functions $f$ on $\hei$ and bounded
$L_\alpha$-harmonic functions $F$ on $S$ by
$$
F(\omega,a)=\int_{\hei}f(\w)P_a^\alpha(\w^{-1}\omega)\
d\w=f*P_a^\alpha(\omega).
$$

We are now in position to prove the following theorem, which gives
characterizations of boundary values of holomorphic or
pluriharmonic $L_\alpha$-harmonic functions in terms of
differential equations. Let us mention that the existence of such
differential equations is well known with different assumptions on
the function (see for instance \cite{L1} or \cite{Gr1,Gr2}), and the
fact that it may be generalized to bounded functions is implicit
in the work of different authors, at least when $\alpha=1$, which
corresponds to the invariant Laplacian. Our aim, here, is to give
a complete proof for bounded functions. It may simplify previous
proofs where such boundary differential
equation appeared.

\begin{theorem} Let $F$ be a bounded $L_\alpha$-harmonic function on
$\uu_n$ with boundary value $f$.
\begin{enumerate}
\renewcommand{\theenumi}{\roman{enumi}}
\item\label{th1:hol} Then $F$ is
holomorphic if and only if  $\lll_nf=0$.
\item\label{th1:anti} Then $F$ is
anti-holomorphic if and only if   $\lll_{-n}f=0$.
\item\label{th1:pluri} Then $F$ is
pluriharmonic if and only if
$\lll_{-n}\lll_nf=(\lll^2+n^2T^2)f=0$.
\end{enumerate}
\label{holph}
\end{theorem}

\begin{proof}  In each case, one implication is elementary. We will content
 ourselves to prove the converse, that is, solutions of the boundary equations
 are holomorphic (resp. anti-holomorphic, resp. pluriharmonic).

 The key point is the following lemma.
 \begin{lemma}
\label{boundaryCR} Let $f$ be a  smooth bounded function on $\hei$
such that $\lll_n f=0$. Then $f$ satisfies the boundary
Cauchy-Riemann equations, that is, \begin{equation} \overline{Z}_k
f=0 \ \ {\mathrm {for}} \ \  k=1,\ldots n.
\label{CR}\end{equation}
\end{lemma}

\begin{proof}
 Replacing eventually $f$ by a left translate of $f$, it is
sufficient to prove it at the point $[0,0]$.

Let us now consider $\ffi$ a smooth compactly supported function
on $\hei$ such that $\ffi=1$ in a neighborhood of $[0,0]$. Let
$\ffi_R(\zeta,t)=\ffi(R^{-1}\zeta,R^{-2}t)$. We may apply Formula
(\ref{Cauchy}) to $\ffi_Rf$ to get
$\ffi_Rf=C(\ffi_Rf)+\widetilde{S}(\lll_n \ffi_Rf)$. The first term
satisfies the boundary Cauchy-Riemann equations, and we are
reduced to consider the second term.  As $\lll_n(\ffi_Rf)=0$ in
some anisotropic ball of radius $cR$, it is equal to some
$\psi_R\lll_n(\ffi_Rf)$ where
$\psi_R(\zeta,t)=\psi(R^{-1}\zeta,R^{-2}t)$ and $\psi=0$ in the
ball of radius $c$. Finally,
\begin{align}
\overline{Z}_kf(0,0)=&\overline{Z}_j\ffi_Rf(0,0)=
\overline{Z}_k\widetilde{S}\bigl(\psi_R\lll_n(\ffi_Rf)\bigr)(0,0)\notag\\
=&\int_{\hei}\overline{Z}_k\Phi(\eta,s)\lll_n(\ffi_Rf)(\eta,s)\psi_R(\eta,s)d\eta
ds\notag\\
=&\int_{\hei}\lll_{-n}(\psi_R\overline{Z}_k\Phi)(\ffi_Rf)d\eta
dt.\notag
\end{align}
But $\ffi_R f$ is bounded, while an elementary computation gives
that $\lll_{-n}(\psi_R\overline{Z}_k\Phi)$ is bounded by
$R^{-3-2n}$. Since we integrate over a shell of volume $R^{2n+2}$,
letting $R\to\infty$, we get $\overline{Z}_k f(0,0)=0$.
\end{proof}

Let us now come back to the proof of the theorem (part
(\ref{th1:hol})). Let us consider a function $F$ such that its
boundary value $f$ is annihilated by $\lll_n$ and show that it
 is holomorphic. We claim that it is sufficient to find a bounded
holomorphic function $G$ that has the same boundary values: indeed, $G$ is
also $L_\alpha$-harmonic, and so $G=F$. Let us now take for $G$
the unique bounded solution of the equation $\lll G+n\partial_a
G=0$ on $\uu_n$, with $G$ given by $f$ on the boundary.
Commutation properties of $\lll$, $T$ and $\partial_a$  give that
$\lll_nG=0$. It follows directly, using the two equations
satisfied by $G$, that $\partial_a G-iTG=0$, that is , $\overline
Z_{n+1}G=0$. Moreover, since on each level set $\{r(z)=a_0\}$ the
function $G$ is smooth, bounded, and satisfies the equation
$\lll_nG=0$ , we also know from the lemma that $\overline{Z}_k G=
0$ for $k=1, 2, \cdots, n$. We have proved that $G$, and so $F$,
is holomorphic.

\medskip

Part (\ref{th1:anti}) is obtained by conjugation.

Let us now turn to assertion (\ref{th1:pluri}) of Theorem
\ref{holph}. We know from part (\ref{th1:anti}) that
$\lll_nf$ is anti-holomorphic. First, note that $\lll_n
Z_{n+1}F=Z_{n+1}\lll_n F=0$. It follows from part (\ref{th1:hol}) that
$$
\overline{Z}_jZ_{n+1}F=0\quad\mbox{for }j=1,\ldots,n+1.
$$
Further, for $k=1,\ldots,n$, note that
\begin{align}
\left[\frac{\partial}{\partial z_k},Z_k\overline{Z}_k\right]F=&
\bigl[(Z_k-2i\overline{z}_kZ_{n+1}),Z_k\overline{Z}_k\bigr]F\notag\\
=&Z_k[Z_k,\overline{Z}_k]F+2i(Z_k\overline{z}_k)(\overline{Z}_kZ_{n+1}F)
+Z_k\bigl((\overline{Z}_k\overline{z}_k)Z_{n+1}\bigr)F.\notag
\end{align}
But $[Z_k,\overline{Z}_k]=-2iT$, $\overline{Z}_kZ_{n+1}f=0$ and
$\overline{Z}_k\overline{z}_k=1$
so that
$$
\left[\frac{\partial}{\partial z_k},Z_k\overline{Z}_k\right]F=
-2iZ_kTF+2iZ_kZ_{n+1}F=-2iZ_kTf+2iZ_k(Z_{n+1}+\overline{Z}_{n+1})F=0.
$$
It follows that $\lll_n\frac{\partial}{\partial
z_k}F=\frac{\partial}{\partial z_k}\lll_nF=0$
which by part (\ref{th1:hol}) implies that $\frac{\partial}{\partial
z_k}F$ is
holomorphic so that $F$ is pluriharmonic.
\end{proof}
\begin{remark}\ This theorem is false if one removes the boundedness property of $F$.
For example, $F[\zeta,t,a]=a^{n\alpha+1}$ is $L_\alpha$-harmonic
on $\hil\times\R$. Its boundary value $f[\zeta,t]=0$
satisfies $(\lll^2+T^2)f=0$ on the boundary but is not pluriharmonic.
\end{remark}

% \begin{remark}\ \\
%--- A much simpler proof is available if $f\in\ss(\hei)$ and not just
%bounded, see \cite{St}
%page 616.\\
%--- This theorem is false if one removes the boundedness property of
%$F$.
%For example, $F[\zeta,t,a]=\abs{\zeta}^2+t-\frac{a}{\alpha}$ is
%$L_\alpha$-harmonic
%on $\hil\times\R$. Its boundary value $f[\zeta,t]=\abs{\zeta}^2+t$
%satisfies $(\lll^2+T^2)f=0$ on the boundary but is not pluriharmonic.
%\end{remark}
%
%In the sequel, our aim is to show that under some mild assumptions the operator
%$L_\alpha$ essentially
%splits into two parts. More precisely, we want to prove that if
%$F$ is a bounded $L_\alpha$-harmonic which satisfies some boundary
%regularity conditions, then $(\lll+n\partial_a)F=0$ and
%$(\partial_a^2+T^2)F=0$. Let us remark that these two conditions
%imply pluri-harmonicity since, from the commutation of the two
%operators $\lll$ and $\partial_{a}$, the first one implies that
%$(\lll^2-n^2\partial_a^2)F=0$. We can state it as a corollary.
%
%\begin{corollary} Let $F$ be bounded function on $\uu_n$ that
%satisfies one of the following
%(equivalent) conditions
%\begin{enumerate}
%\renewcommand{\theenumi}{\roman{enumi}}
%\item $(\lll-n\partial_a)F=(\partial_a^2+T^2)F=0$, or
%\item\label{cor:1:2} $L_\alpha F=L_\beta F=0$ for some
%$\alpha\not=\beta>0$,
%\end{enumerate}
%then $F$ is pluriharmonic.
%\end{corollary}
\subsection{Pluriharmonicity of centrally-independent bounded
boundary-regular functions}

In a first step, we will consider $L_\alpha$-harmonic functions
$f$ that do not depend on the central variable $t$, that is,
functions which satisfy $L_\alpha f=0$ and $Tf=0$. Since $t$ acts
as a parameter, it follows directly from the analysis of the
corresponding situation on $\C^n$. We first compute explicitely
the corresponding Poisson kernel.

\begin{notation} Let $\Delta$ be the standard Laplacian on $\C^n$
(identified with $\R^{2n}$). For $\alpha>0$, let
$$
\Lambda_\alpha=\alpha a(\Delta-n\partial_a)+a^2\partial_a^2
$$
be an operator on $\C^n$. Let $Q_a^\alpha$ be the
Poisson kernel  for $\Lambda_\alpha$. A bounded
$\Lambda_\alpha$-harmonic function $F$ over $\C^n$ is given by
$$
F(\zeta)=f*_{\C^n}Q_a^\alpha(\zeta)=
\int_{\C^n}f(\zeta')Q_a^\alpha(\zeta-\zeta')d\zeta'
$$
with $f\in L^\infty(\C^n)$.
\end{notation}

Let us first identify the Fourier transform over $\C^n$ of
$Q_a^\alpha$.

\begin{lemma}
\label{pre123} The Fourier transform $\widehat{Q_a^\alpha}$ of the kernel $Q_a^\alpha$ is given by
$\widehat{Q_a^\alpha}(\xi)=z(\alpha|\xi|^2a)$, where $z$ is
the unique bounded solution of the equation
\begin{equation}
\bigg(-1-\alpha n\partial_a+a\partial_a^2\bigg)z(a)=0
\label{eqpre123}
\end{equation}
with $z(0)=1$.
In particular, $\widehat{Q_a^\alpha}$ is a smooth function
on
$\C^n\setminus\{0\}$.
\end{lemma}

\begin{proof} The properties of solutions of Equation
(\ref{eqpre123}) are given in Appendix \ref{appendixb}. In
particular, this equation has a unique bounded solution with $z(0)=1$.

Next, observe that $a\to\widehat{Q_a}(\xi)$ is a bounded
solution of the equation
\begin{equation}
\bigg(-\alpha|\xi|^2-\alpha n\partial_a+a\partial_a^2\bigg)y(a)=0.
\label{1.25}
\end{equation}
With the change of variables $y(a)=z(\alpha\abs{\xi}^2a)$, $\xi$
fixed, this equation transforms into (\ref{eqpre123}) so that
there exists $c(\xi)$ such that $\widehat{\widetilde
P_a}(\xi)=c(\xi)z(\alpha|\xi|^2a)$.

Further
$$
c(\xi)=\lim_{a\to 0}\widehat{Q_a}(\xi)=1,
$$
which allows to conclude.
\end{proof}
Let us now state the main result of this section.

\begin{proposition}\label{prop:pluri2} Let $\alpha>0$ and $k$ be the
smallest integer bigger than
$n\alpha$. Let $F$ be a bounded $\Lambda_\alpha$-harmonic function
on $\C^n$ that satisfies the following boundary regularity
condition: for every $p=0,\ldots,k+1$ and every
$\ffi\in\ss(\C^n)$,
$$
\sup\limits_{a\leq1}\abs{\int_{\C^n}\partial_a^p
F(\zeta,a)\ffi(\zeta)d\zeta}<\infty\ .
$$
Then $F$ is constant.
\end{proposition}

\begin{proof} Let us write
 $F=f*_{\C^n}Q_a^\alpha$, with $f$ a bounded function.
Now, let $\ffi\in\ss(\C^n)$ be a function such that
$\widehat{\ffi}$ has compact support with $0\notin\supp
\widehat{\ffi}$. We claim that it is sufficient to prove that, for
all such $\ffi$, we have the identity $
\int_{\C^n}F_a(\zeta)\ffi(\zeta)d\zeta=0$. Indeed, as a
consequence, the Fourier transform of $f$ is supported in $\{0\}$.
It follows that $f$ is a polynomial. Since $f$ is bounded, then
$f$ is constant and finally, $F$ is  constant.

By definition of the Fourier transform of $f$ as a distribution,
we may write
$$
\int_{\C^n}F_a(\zeta)\ffi(\zeta)d\zeta=\scal{\widehat{f},\widehat{\ffi}\widehat{\widetilde
Q^\alpha_a}}
$$
whence,
\begin{equation}
I_p(a):=\partial_a^p\int_{\C^n}F_a(\zeta)\ffi(\zeta)\,d\zeta=
\partial_a^p\scal{\widehat{f},\widehat{\ffi}\widehat{\widetilde
Q_a^\alpha}}
=\alpha^p\scal{\widehat{f}(\cdot),\widehat{\ffi}(\cdot)|\cdot|^{2p}z^{(p)}(|\cdot|^2a)},
\label{1.24}
\end{equation}
where the function $z$ is given in the previous lemma. From now
on, we assume that $n\alpha$ is not an integer. It is then easy to
modify the following proof to cover the remaining case. According to Appendix \ref{appendixb},
for $p=k+1$, we have
$$
\partial_a^{p}z(\alpha|\xi|^2a)=\gamma(|\xi|^2a)+
a^{n\alpha-k}|\xi|^{2(n\alpha-k)}\widetilde{\gamma}(|\xi|^2a),
$$
with $\gamma$, $\widetilde{\gamma}$ smooth functions up to $0$ of
order $N$ ,  and $\widetilde{\gamma}(0)\not=0$. We choose $N$
large enough, depending on the order of the distribution
$\widehat{f}$. Then
$$
I(a)\simeq ca^{n\alpha-k}\,,\ \ \mathrm{with} \ \
c=\widetilde{\gamma}(0)\scal{\widehat{f}(\cdot),\widehat{\ffi}(\cdot)|\cdot|^{2n\alpha+2}},
$$
unless the constant $c$ vanishes.
By assumption, it is a bounded function of $a$. So $c=0$, that is,
$\scal{\widehat{f}(\cdot),\widehat{\ffi}(\cdot)|\cdot|^{2n\alpha+2}}=0$.

Since any function in $\ss(\C^n)$ with compact support in
$\C^n\setminus\{0\}$ can be written as
$\widehat{\ffi}(\cdot)|\cdot|^{2n\alpha+2}$, we conclude that
$\widehat{f}$ vanishes outside $0$, which we wanted to prove.

\end{proof}

\begin{remark} We may identify functions on $S$ that are
independent on  $t$ with functions on the
hyperbolic upper half-plane $\R^{2n}\times\R^+_*$. This one
is an unbounded realization of the real hyperbolic space of
dimension $2n+1$,  while a bounded realization is given by
 the real hyperbolic ball $\B_{2n+1}$ of dimension
$2n+1$. The above theorem is
therefore an analogue
of Theorem 8 in \cite{Ja1} (see also Proposition 2 in \cite{Ja2}),
stating that hyperbolic-harmonic
functions on $\B_{2n+1}$ that are regular up to the boundary are
constant.
\end{remark}

\subsection{Representations of the Heisenberg group and the Poisson
kernel for $L_\alpha$}\ \\
\label{subsec:rephei}
We now pass to the general case. We will have to use harmonic analysis
as before, the scheme of the proof being somewhat similar. Let us
first recall some definitions which are linked to harmonic analysis
on the Heisenberg group.

For $\lambda\neq 0$, the unitary representation $R^\lambda$ of $\hei$ on
$L^2(\R^n)$ is defined by
$$
R^\lambda(\zeta,t)\Phi(x)=e^{2\pi
i\lambda(u.x+u.v/2+t/4)}\Phi(x+v)
$$
if we set $\zeta=u+iv$. The Fourier transform of an integrable
function $f$ on $\hei$ is then the operator-valued function
$\lambda\mapsto\hat{f}(\lambda)$ given by
$$
\scal{\hat
f(\lambda)\phi,\psi}=\int_{\hei}\scal{R^\lambda(\zeta,t)\phi,\psi}f(\zeta,t)d\zeta
dt.
$$

\begin{notation}
Recall that the Hermite functions  $H_k$ of one variable are
defined by the formula
$$
H_k(x)=(-1)^ke^{x^2/2}\left(\frac{d}{dx}\right)^ke^{-x^2},\quad
k=0,1,2,\ldots.
$$
For $c_k=(\sqrt{\pi}2^kk!)^{1/2}$, we note $h_k=c_kH_k$. This
family of functions forms an orthonormal basis of $L^2(\R)$. For
$k=(k_1,\ldots,k_n)$ a multi-index, we write $h_k=h_{k_1}\ldots h_{k_n}$.

Now, for $\lambda\not=0$, denote by $h_k^\lambda(x)=
(2\pi\abs{\lambda})^{n/4}h_k\bigl((2\pi\abs{\lambda})^{1/2}x\bigr)$.
Finally, let
$$
e^{\lambda}_\kappa(\omega)=\sum_{|k|=\kappa}(R^{\lambda}(\omega)
h_k^{\lambda},h_k^{\lambda}).
$$
These last functions may be used to write explicitly an inverse
Fourier Formula (see \cite{T}, or Appendix C).

Let $P_a^\alpha$ be the Poisson kernel for $L_\alpha$, i.e. the
function on $\hei$ that establishes a one-to-one correspondence
between bounded functions $f$ on $\hei$ and bounded
$L_\alpha$-harmonic
functions $F$ on $S$ by
$$
F(\omega,a)=\int_{\hei}f(\w)P_a^\alpha(\w^{-1}\omega)\
d\w=f*P_a^\alpha(\omega),
$$
normalized by $\norm{P_a^\alpha}_{L^1(\hei)}=1$.
Further, for a function $f$ on $\hei$ and $\omega\in\hei$, write
${}_\omega f$ for the function given by
${}_\omega f(\w)=f(\omega\w)$.
\end{notation}

\begin{lemma}
\label{1.4}
Let $\lambda\neq 0$, $\kappa$ be an integer, $\alpha>0$ and $\rho=\alpha(2\kappa+n)$.
Fix $\omega\in\hei$   and define
$$
g_\kappa^\lambda(\omega,a)=\sum_{k\,:\
\abs{k}=\kappa}\scal{\widehat{{}_\omega P_a^\alpha}(\lambda)
h_k^\lambda,h_k^\lambda}=\int_{\hei} {_\omega}
P_a^\alpha(\w)e_\kappa^{\lambda}(\w)\,d\w.
$$
Then
$g^{\lambda}_\kappa(\omega,a)=e^{\lambda}_\kappa(\omega^{-1})g(|{\lambda}|a)$,
where $g$ is the unique bounded solution on $\R^+$ of the equation
\begin{equation}
\label{eq1.4bis}
\Bigg(\partial_a^2-\frac{\alpha n}{a}
-\bigg(\frac{\rho}{a}+1\bigg)\Bigg)g(a)=0
\end{equation}
with $g(0)=1$.
\end{lemma}

\begin{proof}[Proof of Lemma \ref{1.4}]

As $P_a^\alpha\in L^1(\hei)$ with $\norm{P_a^\alpha}_{L^1(\hei)}=1$
and as
$e_\kappa^\lambda$ is easily seen to be bounded, $g^{\lambda}_\kappa$
is bounded.
Further
$$
\lll\,e^{\lambda}_\kappa=(2\kappa+n)|{\lambda}|e^{\lambda }_\kappa.
$$
Since $(\w,a)\to P_a^\alpha(\omega\w)$ is $L_\alpha$-harmonic, we
have
$$
\int_{\hei}
L_\alpha(\,{_\omega}P_a^\alpha)(\w)e^{\lambda}_\kappa(\w)\ d\w=0.
$$
By Harnack's inequality, we may interchange the integral with
the corresponding differential operators to obtain
\begin{equation}
\int_{\hei}\left(\partial_a^2-\frac{\alpha n}{a}\partial_a\right)\,
{_\omega}P_a^\alpha(\w)
e_\kappa^{\lambda}(\w)\
d\w=\left(\partial_a^2-\frac{\alpha
n}{a}\partial_a\right)g^{\lambda}_\kappa(\omega,a),
\label{1.7}
\end{equation}
\begin{equation}
\int_{\hei}\lll\,{_\omega}P_a(\w)e_\kappa^{\lambda}(\w)\
d\w=\int_{\hei}  {_\omega} P_a(\w)\lll e_\kappa^{\lambda}(\w)\
d\w=(2\kappa+n)|{\lambda}|g^{\lambda}_\kappa(\omega,a),
\label{1.8}
\end{equation}
\begin{equation}
\int_{\hei}T^2 {_\omega} P_a(\w)e_{\kappa}^{\lambda}(\w)\ d\w
=\int_{\hei}{_\omega} P_a(\w)T^2e_\kappa^{\lambda}(\w)\ d\w
=-{\lambda}^2g^{\lambda}_\kappa(\omega,a).
\label{1.9}
\end{equation}
Combining (\ref{1.7})-(\ref{1.9}), we see that
the function
$g(a)=g^{\lambda}_k(\omega,a)$ is a bounded solution of the equation
\begin{equation}
\Bigg(\partial^2_a-\frac{\alpha n}{a}\partial_a-
\bigg(\frac{\rho|{\lambda}|}{a}+\lambda^2\bigg)\Bigg)g(a)=0.
\label{1.5}
\end{equation}
One immediately gets that, for $\omega,\lambda$ fixed,
$G_{\omega,\lambda}$
defined by
$$
G_{\omega,\lambda}(a)=g^{\lambda}_\kappa(\omega^{-1},a/|{\lambda}|)
$$
satisfies equation (\ref{eq1.4bis}). The fact that this equation has
only one bounded solution, up to a constant, is well known (see
Appendix \ref{appendixa}, Lemma \ref{hyper}). So
$g^{\lambda}_\kappa(\omega^{-1},a)=c^{\lambda}_\kappa(\omega)g(|{\lambda}|a)$
for
some function $c^{\lambda}_\kappa$.

Moreover, note that $g^\lambda_\kappa$ may be rewritten
\begin{equation}
g^{\lambda}_\kappa(\omega^{-1},a)=e^{\lambda}_k*\check{P_a^\alpha}(\omega),
\label{1.12}
\end{equation}
so that, using  the fact that $P_a^\alpha$ is an approximate
identity, by
letting $a\to0$ in (\ref{1.12}), we get
$c^{\lambda}_\kappa=e^{\lambda}_\kappa$.
\end{proof}
\begin{remark} Even if we will not use it, let us remark that this
allows to write explicitly (see appendix \ref{heifinv})
$$
P_a^\alpha(\omega)=c_n\int_\R\sum_{\kappa\in\N}g_\kappa^\lambda(\omega,a)
|\lambda|^nd\lambda.
$$
\end{remark}
\subsection{An orthogonality property for bounded boundary-regular
functions}\ \\
\begin{notation}
Given $\psi\in C_c^\infty(\R\setminus\{0\})$ we define
$$
e^\psi_\kappa(\w)=\int_\R
e^{\lambda}_\kappa(\w)\psi({\lambda})\,d{\lambda}.
$$
\end{notation}
Of course $e^\psi_k$ is a Schwartz function on $\hei$.
We are now in position to prove the following:

\begin{proposition}\label{1.10}
Let $f\in L^\infty(\hei)$ and let
$F$ be the corresponding bounded $L_\alpha$-harmonic function. Let $k$ be the
smallest integer bigger then $n\alpha$. Assume that for every
$\ffi\in \ss(\hei)$ and every $0\leq p\leq k+1$,
\begin{equation}
\sup_{a\leq1}\abs{\int_\hei\partial_a^p\,F(\w,a)\ffi(\w)\,d\w}<\infty.
\label{1.11}
\end{equation}
Then, for every $\kappa\neq 0$, every $\omega\in\hei$ and every
$\psi\in C_c^\infty(\R\setminus\{0\})$,
\begin{equation}
\label{ortho}
\int_\hei f(\omega\w)e^\psi_\kappa(\w)\ d\w=0.
\end{equation}
\end{proposition}

\begin{proof}[Proof of Proposition \ref{1.10}] Property (\ref{1.11})
is unchanged if one replaces $f$ by
${_\omega}f$ so that it is enough to consider the case $\omega=[0,0]$.

Define
$$
I(a) =\int_\hei F(\w,a)e^{\psi}_\kappa(\w)\,d\w=
\int_\hei f*P_a(\w)e^{\psi}_\kappa(\w)\,d\w
$$
where $f$ is the boundary value of $F$. With (\ref{1.12}) and
Lemma \ref{1.4}, we get
\begin{align}
I(a)=& \int_\hei f(\w)\int_\R
e^{\lambda}_\kappa*\check{P_a}(\w)\psi(\lambda
)\,d{\lambda}\,d\w
\notag\\
=&\int_\hei f(\w)\int_\R
e^{\lambda}_\kappa(\w)g(|{\lambda}|a)\psi({\lambda})
\,d{\lambda}\,d\w.\notag
\end{align}
Using Lemma \ref{nor} in Appendix C, we see that $I$ is a smooth
function for $a>0$ and get
\begin{equation}
\partial^{k+1}\,I(a)
=\int_\hei f(\w)\int_\R
e^{\lambda}_\kappa(\w)|{\lambda}|^{k+1}\partial_a^{k+1}\,
g(|{\lambda}|a)\psi({\lambda})\,d{\lambda}\,d\w
\label{der-p}.
\end{equation}
From now on, we assume that $\alpha n$ is not an integer. It is again
easy to adapt the proof to the other case. Using Lemma \ref{hyper}
in Appendix A with $N=1$,  can write
$$
\partial^{k+1}\,g(|{\lambda}|a)=|{\lambda}|^{\alpha n-k}a^{\alpha
n-k}g_1(|{\lambda}|a)+ g_2(|{\lambda}|a)
$$
with $g_1$ and $g_2$
having continuous derivatives up to $0$, and $g_1(0)\neq 0$. Using
again Lemma \ref{nor}, we can pass to the limit when $a$ tends to
$0$ in the two integrals. We find that
$$
\partial^{k+1}\,I(a)\simeq c a^{\alpha
n-k} \ \ \mathrm{with} \    \ c=g_1(0)\int_\hei f(\w)\int_\R
e^{\lambda}_\kappa(\w)|{\lambda}|^{k+1}\psi({\lambda})\,d{\lambda}\,d\w
$$
unless $c=0$. Now, Hypothesis (\ref{1.11}) implies that
$\partial^{k+1}_aI(a)$ is bounded, so it is indeed the case and
$$
\int_\hei f(\w)\int_\R
e^{\lambda}_\kappa(\w)|{\lambda}|^{k+1}\psi({\lambda})\,d{\lambda}\,d\w
=0.
$$
Now every function in $\cc^\infty_c(\R\setminus\{0\})$ can be
written in the form  $|{\lambda}|^{k+1}\psi({\lambda})$. This
completes the proof of the proposition.
\end{proof}

\subsection{The main theorem on the Heisenberg group}

We are now in position to  prove our main theorem in this context:

\begin{theorem} Let $\alpha>0$ and $k$ the smallest integer bigger
than $n\alpha$. Let $F$ be a bounded
$L_\alpha$-harmonic function on $S$. Assume that for every
$\ffi\in \ss(\hei)$ and every $0\leq p\leq k+1$,
\begin{equation}
\sup_{a\leq1}\abs{\int_\hei\partial_a^p\,F(\w,a)\ffi(\w)\,d\w}<\infty.
\tag{\ref{1.11}}
\end{equation}
Then, $F$ is pluriharmonic.
\label{1.2}
\end{theorem}

\begin{proof}
Let $f\in L^\infty(\hei)$ be the  boundary value of $F$. We will
use the following well-known fact: {\sl the pointwise limit of a
uniformly bounded sequence of pluriharmonic functions is again
pluriharmonic.} It allows to replace $f$ by its right convolution
with a smooth approximate identity. Moreover, such a function
satisfies also (\ref{1.11}). So, from now on in this proof, we may
assume that $f$ is smooth, and that, moreover, its derivatives up
to order $4$ are bounded.

To prove that $F$ is pluriharmonic, it is sufficient to prove that
$g(\w)=(\lll^2+n^2T^2)f(\w)$ is constant. Indeed, by Harnack's Inequality,
such a function tends to $0$ when $\w$ tends to $\infty$, so $g$
can only be equal to the constant $0$. We then use Theorem
\ref{holph} to conclude.

Let $G$ be the  $L_\alpha$-harmonic extension of $g$. Then $G$
satisfies also hypothesis (\ref{1.11}), and according to
 (\ref{orth}) in Proposition \ref{1.10}, for every $\kappa\not=0$,
every $\omega\in\hei$
and every $\psi\in\cc^\infty_0(\R\setminus\{0\})$ we have
\begin{equation}\label{orth}
\int_\hei g(\w)e_\kappa^\psi(\w)\,d\w=0.
\end{equation}
Moreover, a direct integration by parts, using the fact that
$$
(\lll^2+n^2T^2)e_0^\lambda= 0,
$$
allows to conclude that (\ref{orth}) is also valid for $\kappa=0$.
We deduce from (\ref{FIhei2}) that  the Fourier transform of $g$
in the $t$-variable (in the distributional sense) is supported by
$0$. It follows that $g$ is a polynomial in the $t$-variable and,
as $g$ is bounded, this implies that $g$ is independent of the
central variable $t$. According to Proposition \ref{prop:pluri2},
it follows that $G$ is  constant, which we wanted to prove.
\end{proof}

\subsection{Optimality of the result}

We will now prove a converse of Theorem \ref{1.2}, showing that
the index of regularity given there is optimal. We will actually
prove that, in some weak sense, all $L_\alpha$-harmonic functions
have the regularity just below the limitation given by that
theorem.

Recall that a function $F$ on $S$ is said to have distributional
boundary value if, for every
$\psi\in\cc^\infty_c(\hei)$, the limit
\begin{equation}
\label{bvd}
\lim\limits_{a\to 0}\int_{\hei}F(\w,a)\psi(\w)\,d\w
\end{equation}
exists. Note that, if $F$ has a boundary distribution, so do $\lll F$ and
$TF$.

We may now prove the following:
\begin{theorem}
Let $\alpha>0$ and let $k$ be the smallest integer greater then
$n\alpha$. Assume $F$ is a $L_\alpha$-harmonic function on $S$
with a boundary distribution. Then, for every $p\leq k$,
$\partial^p_a F$ has a boundary distribution. \label{th:bd}
\end{theorem}

\begin{proof}  The proof is essentially
the same as in \cite{Ja2} or \cite{BBG} which deal with bounded
realizations. So we  only give a quick outline of it.

We prove the theorem by induction on $p$. For $p=0$, this is our
original assumption. Let us take the statement for granted at rank
$p-1$ and fix $\psi\in\cc^\infty_c(\hei)$. Define
$$
\psi_p(a)=\partial_a^p\int_\hei F(\w,a)\psi(\w)\,d\w=
\int_\hei\partial_a^p F(\w,a)\psi(\w)\,d\w.
$$
Since
$$
a\partial_a^2F-n\alpha\partial_aF=(\alpha\lll+aT^2)F,
$$
applying $\partial_a^{p-1}$ to both sides, we get
\begin{equation}
\label{th:bd:eq1}
a\partial_a^{p+1}F+(p-1-n\alpha)\partial_a^pF=\partial_a^{p-1}(\alpha\lll+aT^2)F.
\end{equation}
We know that $(\alpha\lll+aT^2)F$ has a boundary distribution. It
follows from the induction hypothesis that the second member of
(\ref{th:bd:eq1}) has a boundary distribution. Now, multiply
(\ref{th:bd:eq1}) by $\psi$ and integrate over $\hei$, we get that
$$
g_p(a):=a\partial\psi_p(a)+(p-1-n\alpha)\psi_p(a)
$$
has a limit as $a\to0$. As $p-1-n\alpha\not=0$, solving this differential equation, we get
$$
\psi_p(a)=\lambda a^{n\alpha-p+1}+a^{n\alpha-p+1}\int_1^a
\frac{g_p(t)}{t^{(n\alpha-p+1)+1}}dt.
$$
As $g_p$ has a limit when $a\to0$, it follows that $\psi_p$ has a
limit as $a\to0$, provided $n\alpha-p+1>0$. We have  proved the
theorem.
\end{proof}

\section{The main theorem on irreducible symmetric Siegel domains of tube type}
\label{sectionmainth}\setcounter{equation}{0}

We first write the Hua operator in an appropriate coordinate
system.

\subsection{Preliminaries on irreducible symmetric cones}
\label{cones}

Let $\Omega$ be an irreducible symmetric cone in an Euclidean
space $V$, as in the introduction. We describe precisely the
solvable group that acts simply transitively on $\Omega$ in terms
of Jordan algebras. We refer to the book of Faraut and Kor\'anyi
\cite{FK} for these prerequisites, or to \cite{BBDHPT} where a
rapid introduction has already been given, with the same
notations.

We assume that $V$, endowed with the scalar product $\langle\cdot
,\cdot \rangle$ is an Euclidean Jordan algebra, that is, is also
endowed with a product such that, for all elements $x,y$ and $z$
in $V$
$$
xy=yx, \hspace{2cm} x(x^2y)=x^2(xy), \hspace{2cm}
\scal{xy,z}=\scal{y,xz}.
$$
Moreover, we assume that $V$ is a simple Jordan algebra with
unit element $e$. We denote by $L(x)$ the self-adjoint
endomorphism of $V$ given by the multiplication by $x$, i.e.
\begin{equation}
L(x)y=xy.
\label{def:L}
\end{equation}

The irreducible symmetric cone $\Omega$ is then given by
$$
\Omega=\,\mbox{int}\,\{x^2:\ x\in V\}.
$$
Let $G$ be the connected component of the group  of all
transformations in
$GL(V)$ which leave $\Omega$ invariant, and let $\mathcal{G}$ be its
Lie algebra. Then $\mathcal{G}$ is a subspace of the space of
endomorphisms of $V$ which contains all
$L(x)$ for all $x\in V$, as well as all $x\,\Box\,y$ for $x,y \in V$,
where
\begin{equation}
\label{def:box}
x\,\Box\,y=L(xy)+[L(x),L(y)]
\end{equation}
(see \cite{FK} for  these properties).

\smallskip

We fix a Jordan frame $\{c_1,\dots,c_r\}$ in $V$, that is, a
complete system of orthogonal primitive idempotents:
$$
c^2_i=c_i,\hspace {2cm} c_ic_j=0 \quad\mbox{if}\; i\neq j, \hspace
{2cm}c_1+...+c_r=e
$$
and none of the $c_1,...,c_r$ is a sum of two nonzero idempotents.
Let us recall that the length $r$ is the {\sl rank} of the cone and is
independent of the choice of the Jordan frame.

The Peirce decomposition of $V$ related to the Jordan frame
$\{c_1,\dots,c_r\}$ (\cite {FK}, Theorem IV.2.1) may be written as
\begin{equation}
\label{Peirce}
V=\bigoplus_{1\leq i\leq j\leq r}V_{ij}\,.
\end{equation}
It is given by the common diagonalization of the self-adjoint
endomorphims $L(c_{j})$ with respect
to their only eigenvalues $0$, $\frac 12$, $1$. In particular
$V_{jj}=\R c_j$ is the
eigenspace of $L(c_{j})$ related to $1$, and, for $i<j$, $V_{ij}$ is
the intersection of the
eigenspaces of $L(c_{i})$ and $L(c_{j})$ related to $\frac{1}{2}$.
All $V_{ij}$, for $i<j$, have
the same dimension $d$.

For each $i<j$, we fix once for all an orthonormal basis of $V_{ij}$,
which we note
$\{e^\alpha_{ij}\}$, with $1\leq \alpha\leq d$. To simplify the
notation, we write
$e^\alpha_{ii}=c_{i}$ ($\alpha$ taking only the value $1$). Then the
system  $\{e^\alpha_{ij}\}$,
for $i\leq j$ and $1\leq \alpha \leq \dim V_{ij}$, is an orthonormal
basis of $V$.

Let us denote by $\aa$ the Abelian subalgebra of $\mathcal{G}$
consisting of elements $H=L(a)$, where
$$
a=\sum_{j=1}^ra_jc_j\in \bigoplus_iV_{ii}.
$$
We set $\lambda_j$ the linear form on $\aa$ given by
$\lambda_j(H)=a_j$. The Peirce decomposition gives also a
simultaneous diagonalization of all $H\in\aa$, namely
\begin{equation}
Hx=L(a)x=\frac{\lambda_{i}(H)+\lambda_{j}(H)}{2}x \hspace{3cm}  x\in
V_{ij}\,.
\label{diag}
\end{equation}
Let $A=\exp\aa$. Then $A$ is an Abelian group, which is the
Abelian part of the Iwasawa decomposition of $G$. We  now describe
the nilpotent part $N_{0}$. Its Lie algebra $\nn_{0}$ is the space
of elements $X\in\mathcal{G}$ such that, for all $i\leq j$,
\begin{equation}
\label{triang}
XV_{ij}\subset \bigoplus_{k\geq l\;;\;( k,l)>(i,j)}V_{kl},
\end{equation}
where the pairs are ordered lexicographically. Once $\nn_{0}$ is defined,
we define $\ss_{0}$
as the direct sum $\nn_{0}\oplus \aa$. The groups $S_{0}$ and $N_{0}$
are then
obtained by taking the exponentials. It follows from the definition
of $\nn_{0}$ that the
matrices of elements of $\ss_{0}$ and $S_{0}$, in the orthonormal
basis
$\{e^\alpha_{ij}\}$, are upper-triangular.

The solvable group $S_{0}$ acts simply transitively on $\Omega$. This
may be found in \cite{FK}
Chapter VI, as well as the precise description of $\nn_{0}$ which
will be needed later. One has
$$
{\nn}_0=\bigoplus_{i<j\leq r}{\nn}_{ij},
$$
where
$$
\nn_{ij}=\{z\,\Box\,c_i\ :z \ \in V_{ij}\}.
$$
This decomposition corresponds to a diagonalization of the adjoint
action of ${\aa}$ since
\begin{equation}
[H,X]=\frac{\lambda_j(H)-\lambda_i(H)}{2}X,\ X\in \nn_{ij}.
\label{comN}
\end{equation}
Finally, let $V^{\C}= V+iV$ be the complexification of $V$. The
action of $G$ is extended to $V^{\C}$ in the obvious way.

\subsection {Preliminaries on irreducible symmetric Siegel domains of tube type}
We keep notations of the previous section and let
\begin{equation}
\dd=\{z\in V^{\C}:\ \Im z\in\Omega\}.
\label{Siegel}
\end{equation}

The elements $x\in V$ and $s\in S_0$ act on
${\dd }$ in the following way:
$$
x\cdot z=z+x\quad,\quad s\cdot z=sz.
$$
Both actions generate  a solvable Lie group
$$
S=VS_0=VN_0A=NA,
$$
which identifies with a group of holomorphic automorphisms acting
simply transitively on  $\dd$.
The Lie algebra $\ss$ of $S$ admits the decomposition
\begin{equation}
\ss=V\oplus\ss_0
=\left( \bigoplus_{i\le j}\,V_{ij} \right)
\oplus \left( \bigoplus_{i<j}\,{\nn }_{ij} \right)
\oplus\aa .\label{direct}
\end{equation}
Moreover, by (\ref{diag}) and (\ref {comN}), one knows the adjoint
action of elements $H\in\aa$:
\begin{align}
[H,X]&=\frac{\lambda_i(H)+\lambda_j(H)}{2}X\quad\text{for}\quad
X\in V_{ij},\label{action}\\
[H,X]&=\frac{\lambda_j(H)-\lambda_i(H)}{2}X\quad\text{for}\quad
X\in{\nn }_{ij}.\nonumber
\end{align}

Since $S$ acts simply transitively on the domain $\dd$,
we may identify $S$ and $\dd$. More precisely, we define
\begin{equation}
\theta: S\ni s\mapsto \theta (s)=s\cdot\mathbf{e}\in \dd
\label{theta},
\end{equation}
where $\mathbf{e}$ is the point $(0,ie)$ in ${\ss}$. The Lie algebra
$\ss$ is then identified with
the tangent space of $\dd$ at $\mathbf{e}$ using the differential
$d\theta_e$. We identify
$\mathbf{e}$ with the unit element of $S$. We transport both the
Bergman metric $g$ and the complex
structure $\jj $ from $\dd $ to $S$, where they become left-invariant
tensor fields on $S$.
We still write $\jj$ for the complex structure on $S$. Moreover, the
complexified tangent
space $T_\mathbf{e}^{\C}$ is identified with the complexification of
$\ss$,
which we denote by $\ss^\C$. The decomposition
$T_\mathbf{e}^{\C}=T_\mathbf{e}^{(1,0)}\oplus T_\mathbf{e}^{(0,1)}$
is transported into
$$
\ss ^\C=\qq\oplus\pp.
$$
Elements of $\ss ^\C$ are identified with left invariant vector
fields on $S$, and are called
left invariant holomorphic vector fields when they belong to $\qq$.
The K\"ahlerian metric
given by the Bergman metric can be seen as a Hermitian form $(\cdot
,\cdot )$
on $\qq$, and orthonormality for left invariant  holomorphic vector
fields means orthonormality
for the corresponding elements in $\qq$.

\medskip
Now, we construct a suitable orthonormal
basis of $\qq $. Let $\{e^\alpha_{jk}\}$ be the orthonormal basis of
$V$
fixed in the previous subsection. For $j<k$ and $1\leq \alpha \leq
d$, we define
$X^\alpha_{jk}\in V_{jk}$ and $Y^\alpha_{jk}\in\nn_{jk}$
as the left-invariant vector fields on $S$ corresponding
to $e^\alpha_{jk}$ and $2e_{jk}^\alpha\,\Box\,c_j$, respectively.
For each $j$ we define $X_j$ and $H_j$ as the left-invariant vector
fields on $S$ corresponding to
$c_{j}\in V_{jj}$ and $L(c_j)\in\aa$, respectively.

Finally, let
$$
Z_{j}=X_{j}-iH_j\quad,\quad Z_{jk}^\alpha=X_{jk}^\alpha
-iY_{jk}^\alpha,
$$
which means that
\begin{equation}
\label{cs}
\jj (X_{j})=H_j\quad,\quad
\jj (X_{jk}^\alpha )=Y_{jk}^\alpha.
\end{equation}

The left invariant vector fields $Z_{j}$, for $j=1,\cdots,r$,
and $Z_{jk}^\alpha$, for $j< k\leq r$  and $\alpha =1,\cdots,d$.
constitute an orthonormal basis of $\qq$.

Using $Z_j$ we can compute the so called {\it strongly diagonal}
Hua operators i.e the operators defined by
$$
\hil _j f=(\hil f\cdot Z_j, Z_j)\quad,\quad j=1,\cdots,r.
$$
In terms of the basis $X_j,X_{jk}^\alpha,Y_{jk}^\alpha,H_j$ they are
\cite{BBDHPT}:
\begin{equation}
\hil_j=\Delta_j+\frac{1}{2}\sum_{k<j}\sum_\alpha
\Delta^\alpha_{kj}+\frac{1}{2}\sum_{l>j}\sum_{\alpha}\Delta^\alpha_{jl},\label{Hj}
\end{equation}
where
\begin{equation}
\label{ceg}
\Delta_j=X_{j}^2+H_j^2-H_j\quad\quad
\Delta_{ij}^\alpha=(X_{ij}^{\alpha })^2+(Y_{ij}^{\alpha })^2-H_j.
\end{equation}
$\Delta_j$'s and  $\Delta_{ij}^\a$'s are $S$-invariant operators, which
at the point $ie$ agree with $\partial_{z_j}\partial_{\overline z_j}$
and $\partial_{z_{ij}^\a}\partial_{\overline z_{ij}^\a}$.

%    \begin{definition}
 %   A function $F$ on $\dd$ is said to be \emph{strongly Hua-harmonic} if
  %  $\hil_1F=\ldots =\hil_rF=0$.
   % \end{definition}
    %Let us remark that this notion is a priori weaker than
    %Hua-harmonicity.

\subsection{Further notations}

In this subsection, we collect some information and some notations
which will be used in the proof, in the next section. We assume
that $r\geq2$. We will define some sub-algebras and subgroups.
Let\arraycolsep 1pt
$$
\begin{array}{ccc}
\aa^-=\mbox{lin}\,\{L(c_1),\ldots,L(c_{r-1})\},&&
\aa^+=\mbox{lin}\,\{L(c_r)\},\\
\nn_0^-=\bigoplus\limits_{i<j\leq r-1}\nn_{ij}&\quad\mbox{and}\quad&
\nn_0^+=\bigoplus\limits_{j=1}^{r-1}\nn_{jr}.\\
\end{array}
$$
Then $\nn_0^+$ is an ideal of $\nn_0$, while $\nn_0^-$ is a
subalgebra.
Clearly
$$
\aa=\aa^-\oplus\aa^+
$$
and
$$
\nn_0=\nn_0^-\oplus\nn_0^+.
$$
Next, we define $A^+$, $A^-$, $N_0^+$, $N_0^-$ as the exponentials of
the corresponding Lie algebras. We have
$$
A=A^-A^+\quad\mbox{and}\quad N_0=N_0^-N_0^+
$$
in the sense that the mappings
$$
\begin{matrix}
A^-\times A^+&\to&A\\
(a^-,a^+)&\mapsto&a^-a^+\\
\end{matrix}
\quad\mbox{and}\quad
\begin{matrix}
N_0^-\times N_0^+&\to&N_0\\
(y^-,y^+)&\mapsto&y^-y^+\\
\end{matrix}
$$
are diffeomorphisms.

%            Let
%            $$
%            V^-=\bigoplus\limits_{1\leq i\leq j\leq r-1}V_{ij}
%            \quad\mbox{and}\quad
 %           V^+=\bigoplus\limits_{j=1}^rV_{jr}.
  %          $$
  %          Then $V^-$ is a subalgebra of $V$ with the Jordan frame
   %         $c_1,\ldots,c_{r-1}$. Let
    %        $$
    %        \Omega^-=\mbox{int}\,\{x^2\,:\ x\in V^-\}
     %       $$
      %      be the corresponding symmetric cone of rank $r-1$. Then clearly
       %     $$
       %     S_0^-=N_0^-A_0^-
        %    $$
         %   is the solvable group corresponding to $\Omega^-$ described in
          %  Subsection \ref{cones}.
           %Indeed, the transformations $L(c_j)|_{V^-}$,
            %$j=1,\ldots,r-1$ and $z\,\Box\,c_j|_{V^-}$, $z\in V_{jk}$,
            %$k=2,\ldots,r-1$, $j<k$ are those define by
            %(\ref{def:L}-\ref{def:box}) for the smaller cone.
 %           The Siegel tube domain of rank $r-1$ is
  %          $$
   %         \dd^-=V^-+i\Omega^-.
    %        $$
    %        The group
     %       $$
      %      S^-=V^-S_0^-
       %     $$
        %    acts simply transitively on $\dd^-$ and it is identified with it by
         %   $$
 %           \begin{array}{cccc}
  %          \theta^-\,:&S^-&\to&\dd^-\\ &s^-&\mapsto&is^-\circ\mathbf{e}^-\\
  %          \end{array}
   %         $$
    %        where $\mathbf{e}^-=\bigl(0,i(c_1+\ldots+c_{r-1})\bigr)\in\dd^-$.
%
 %           In addition, if
  %          $$
   %         S^+=V^+N_0^+A^+,
    %        $$
     %       then
      %      $$
       %     S=S^+S^-=S^-S^+
       %     $$
        %    in the same sense as before.

%$$
%\begin{matrix}
%S^+\times S^-&\to&S\\
%(s^+,s^-)&\mapsto&s^+s^-\\
%\end{matrix}
%\quad\mbox{and}\quad
%\begin{matrix}
%S^-\times S^+&\to&S\\
%(s^-,s^+)&\mapsto&s^-s^+\\
%\end{matrix}
%$$
%are diffeomorphisms.

\subsection{Special coordinates}
Let   $\partial\dd = V + i\partial\Omega$ be the topological
boundary of the domain $\dd = V + i\Omega$. We fix a Jordan frame,
choose coordinates in $V$ according to the Peirce decomposition $x
= \sum x_{jk}^\a e_{jk}^\a$ and we order them lexicographically
i.e. $(j,k)>(l,p)$ if either $j>l$ or $j=l$ and $k>p$.

Let $S=VN_0A$ be the corresponding solvable Lie group. We consider its subgroup
$S' = VN_0A^-$ that will be identified with $\R^{2n-1}$.
To determine appropriate coordinates in $S'$ we look more carefully at the transformations
that build $N_0$. Namely \cite{FK}, any element of $N_0$ can be written uniquely as
$$\tau(y^1)\cdot\ldots\cdot\tau(y^{r-1}),$$
where $y^{j}= \sum_{\substack{j < k \le r\\ \a}}y^\a_{jk}e^\a_{jk}$
and $\tau(y^j)=\exp(2y^j\Box c_j)$.

Therefore, any element
$x\tau(y^1)\cdot\ldots\cdot\tau(y^{r-1})\exp(\sum_{j < r}y_{jj}L(c_j))$ of $S'$
can be uniquely written as
$(x,y_{11},y^1,y_{22},y^2,\ldots,y^{r-1}).$ Notice that the
coordinates are ordered lexicographically
exactly as  in $\dd\subset V+iV$. Now we define a diffeomorphism $\phi$ of $S'=\R^{2n-1}$
onto an open subset of $\partial\dd$ as the limit point of the curve $t\mapsto s'\exp(- tH_r)\cdot ie$:
$$
\phi(s')=\lim_{t\to\infty}s'\exp(- tH_r)\cdot ie.
$$
In the above coordinates we have
\begin{align}
\phi(x,y_{11},y^1,y_{22},y^2&,\ldots,y_{r-1,r-1},y^{r-1})
=\notag\\
&=x + i \left(\sum_{\substack{j< k\\  \alpha}}
\big(e^{y_{jj}}y_{jk}^\a + P_{jk}^\a(y)\big)e_{jk}^\alpha
+i\sum_{j < r}\big(e^{y_{jj}}+P_{jj}(y)\big)c_j+ P_{rr}(y)c_r
\right),
\label{lgeod}
\end{align}
where $P_{jk}^\a$ and $P_{jj}$ are polynomials depending only on
the coordinates that proceed $y_{jk}^\a$ or $y_{jj}$ in the above
order i.e. on $y_{lp}$ for $(l,p)<(j,k)$ or $(l,p)<(j,j)$, which
follows from triangularity of the action of the group $N_0$, given
by  \eqref{triang}. To obtain explicit formulas for the
polynomials one has to use a more precise formula for $\tau$, as
described in Chapter VI of \cite{FK}, see also \cite{DHMP}.

\begin{lemma}
The mapping $\phi$ is an one-to-one diffeormorphism of $\R^{2n-1}$
onto an open subset of $\partial\dd$.
\end{lemma}

\begin{proof} One gets that
$\phi$ is one-to-one from formula (\ref{lgeod}) and the observation that the action of $N$
is triangular. To prove that $\phi$ is a diffeomorphism it is enough to compute $d\phi$,
the differential of $\phi$.
Using again triangularity of the action of $N$,
one can easily find a minor
of rank $2n-1$, which is a triangular matrix with $1$'s and $e^{y_{jj}}$'s
on the diagonal. The Inverse Mapping Theorem
implies the claim.
\end{proof}

Using $\phi$ we define a coordinate system on a neighborhood of $\phi(\R^{2n-1})$
$$
\R^{2n-1}\times\R\ni(w,b)\mapsto\phi(w)+ibc_r\in V+iV,
$$
where $w=(x,y_{11},y^1,y_{22},y^2,\ldots,y_{r-1,r-1},y^{r-1})$ is identified with
the corresponding element $s'$ of the group $S'$.
%$$x\tau(y^1)\circ\ldots\circ\tau(y^{r-1})\exp(\sum_{1\le j< r}y_{jj}L(c_j))\in S'.$$
This
means that for positive $b$
\begin{equation}
\label{star} \phi(w)+ibc_r=s'\exp (bL(c_r))\cdot ie\in\dd
\end{equation}
and
\begin{equation}
\label{2star} \phi(w)+ibc_r\notin\overline\dd
\end{equation}
if $b$ is negative. For every Jordan frame $c_1,\ldots,c_r$ and
the corresponding group $S$ we may construct such a system.
Moreover, applying an element $g$ of the group $G$ to it, we
obtain a coordinate system on a neighborhood of $g{(\phi(S'))}$
satisfying \eqref{star} and  \eqref{2star}. We are going to
exploit such systems to define regularity of a function near the
boundary of $\dd$.

\begin{definition}
We say that a coordinate system $\Phi:\R^{2n-1}\times\R\mapsto
V+iV$ is \emph{ a special coordinate system} if it is of the form
\begin{equation}
\label{o}
\Phi(w,b)=g\big(\phi(w)\big)+ibg(c_r)
\end{equation}
for a Jordan frame $c_1,...c_r$ and a $g\in G$.
\end{definition}

Special coordinate systems are suitable to describe the boundary
behavior of bounded pluriharmonic (holomorphic) functions, in
terms of some integral conditions. More precisely, we consider
functions satisfying the following regularity condition:
\vskip3pt
\noindent{\bf Condition \ref{cond1}.} {\sl A function $F$ is said to satisfy
Condition \ref{cond1} for some integer $k$ if, for every special
coordinate system $\Phi$, and every $\psi\in\ss (\R^{2n-1})$,
\begin{equation}
\sup_{0<b<1}\abs{\int_{\R^{2n-1}}\partial^k_{b}F\big(\Phi(w,b)\big)\psi(w)\,
dw}<\infty.
\label{cond1}
\end{equation}
for some integer number k, where $dw$ is the Lebesgue measure.}

\begin{proposition}
Let $F$ be a bounded pluriharmonic  function on $\dd$, then $F$
satisfies Condition \ref{cond1} for every $k$.
\end{proposition}

\begin{proof}
The action of $G$ preserves pluriharmonicity and $F(\Phi( w,b)) =
(F\circ g)(\phi(w)+ibc_r)$ for some $g\in G$. Notice also that if
$F$ satisfies Condition \ref{cond1} for some $k$, then so does
$F_g(z)=F(g\cdot z)$. Therefore, we may assume that $g=Id$ in
\eqref{o}.
Moreover, it is enough to prove \eqref{cond1} for even $k$'s.\\
\medskip
Now since
$$
0 = \partial_{z_r}\partial_{\overline z_r}F =
 (\partial_{x_r}^2+\partial_{b}^2)F,
$$
we have
\begin{eqnarray*}
\int_{\R^{2n-1}}\partial^{2k}_{b}F(\phi(w)+ibc_r)\psi(w)dw
&=& \int_{\R^{2n-1}} (-\partial^{2}_{x_r})^k F(\phi(w)+ibc_r)\psi(w)dw\\
&=& (-1)^k\int_{\R^{2n-1}} F(\phi(w)+ibc_r)(\partial^{2}_{x_r})^k \psi(w)dw
\end{eqnarray*}
and \eqref{cond1} follows.
\end{proof}

\subsection{Maximum regularity of Hua-harmonic functions}

Now we are in position to formulate and to prove the main result of the paper

\begin{mth}
\label{mainth}
Let $\Omega$ be an irreducible symmetric cone if rank $r$
and let $d$ be the common dimension of the associated spaces in its Peirce
decomposition (\ref{Peirce}). Let
$\dd$ be the tube type domain associated to $\Omega$
If $F$ is a bounded Hua-harmonic function on $\dd$
satisfying Condition \ref{cond1} for
$k=0,\ldots,\ent{\frac{(r-1)d+1}{2}}+1$, then $F$ is
pluriharmonic.
\end{mth}
\begin{remark}
Notice that $F$ in the above theorem has certain boundary
regularity at any point $z\in\partial\dd$ for which $\Im z$
has rank $r-1$ in the Jordan algebra. The set of such points will be
denoted $\widetilde{\partial\dd}$. The theorem says that this
regularity implies pluriharmonicity.
\end{remark}
Fix a Jordan frame $c_1,\ldots,c_r$, the corresponding group $S$ and
differential operators $\Delta_k, \Delta_{jk}^\alpha$, \eqref{ceg}.
We will show that\arraycolsep 1.5pt
\begin{equation}
\begin{array}{ccccl}
\Delta_kF&=&0&\quad&\mbox{for }k=1,\ldots,r\\
\Delta_{jk}^\alpha F&=&0&&\mbox{for }1\leq j<k\ldots,r\mbox{ and }
\alpha=1,\ldots,d\\
\end{array}
\label{cond2}
\end{equation}
and the conclusion will follow by \cite{BDH}.

\medskip
Let
$$
N^+=V^+N_0^+,\quad N^-=V^-N_0^-.
$$
Then clearly,
$$
N=N^-N^+,\quad S^-=N^-A^-,\quad S^+=N^+A^+.
$$
Denote by $dn^-$, $dn^+$, $da^-$, $da^+$ the Haar measures on the
groups
$N^-$, $N^+$, $A^-$, $A^+$ respectively and $ds^-=dn^-da^-$,
$ds^+=dn^+da^+$.

Notice that $N^+$ is the Heisenberg group and $S^+$ is one parameter
extension considered in Section \ref{heisenberg}. Indeed, for the
basis
$X_{jr}^\alpha$, $Y_{jr}^\alpha$, $X_r$, $H_r$ of the Lie algebra of
$S^+$ we have
$$
\begin{array}{rclcrcl}
\bigl[Y_{jr}^\alpha,X_{jr}^\alpha\bigr]&=&X_r,&\qquad&
\bigl[H_r,X_r\bigr]&=&X_r,\\[3mm]
\bigl[H_r,X_{jr}^\alpha\bigr]&=&\frac{1}{2}X_{jr}^\alpha,&&
\bigl[H_r,Y_{jr}^\alpha\bigr]&=&\frac{1}{2}Y_{jr}^\alpha,\\
\end{array}
$$
while all other brackets are $0$ (see (1.16) and (1.17) of
\cite{DHMP}).

Given $\psi\in\cc_c^\infty(S^-)$, let
$$
G_\psi(s^+)=\int_{S^-}F(s^-s^+)\psi(s^-)\,ds^-.
$$
The operator $\hil_r$ is well defined for a function on $S^+$ and
\begin{equation}
\hil_rG_\psi(s^+)=\int_{S^-}(\hil_rF)(s^-s^+)\psi(s^-)\,ds^-.
\label{eq2.3}
\end{equation}
Therefore
$$
\hil_rG_\psi=0.
$$
Notice that $\hil_r$ is the operator $L_{1/2}$ from Section
\ref{heisenberg}.

\begin{lemma} Assume that $F$ satisfies the conditions of Theorem
\ref{mainth}. Then, for every $\phi\in\ss (N^+)$ and every
$k=1,\ldots,\ent{\frac{(r-1)d+1}{2}}+1$,
\begin{equation}
\sup_{0<a^+<1}\abs{\int_{N^+}\partial_{a^+}^k G_\psi(n^+a^+)
\phi(n^+)\,dn^+}<\infty.
\label{cond4}
\end{equation}
\label{lem2.1}
\end{lemma}

\begin{proof} Define
$$
I(a^+)=\int_{N^+} G_\psi(n^+a^+)\phi(n^+)\,dn^+
$$
so that
\begin{align}
I(a^+)=&\int_{N^+}\int_{S^-}F(s^-n^+a^+)\psi(s^-)\phi(n^+)
\,ds^-dn^+\notag\\
=&\int_{N^+}\int_{S^-}F(n^-a^-n^+a^+)\psi(s^-)\phi(n^+)
\,dn^-da^-dn^+.\notag
\end{align}
Now we change variables via the transformation
$$
n^+\mapsto(a^-)^{-1}n^+a^-.
$$
Then
\begin{eqnarray*}
I(a^+)&=&\int_{S^-}\int_{N^+}F(n^-n^+a^-a^+)\psi(s^-)\phi
((a^-)^{-1}n^+a^-)
\prod _{j=1}^{r-1}a_j^{-\frac{d}{2}}\,dn^+dn^-da^-\\
&=& \int_{\R^{2n-1}}F(\Phi(w,a^+))\widetilde\psi(w)dw
\end{eqnarray*}
for some $\widetilde\psi$ in the Schwartz class. So (\ref{cond4})
follows from (\ref{cond1}).
\end{proof}

\begin{lemma}\label{lem2.2}
Assume that $F$ satisfies the conditions of Theorem \ref{mainth}, then
\begin{equation}
\label{delr}
\Delta_rF=0\quad\text{and}\quad
\Delta_{jr}^\alpha F=0
\end{equation}
for $j=1,\ldots,r-1$ and $\alpha=1,\ldots,d$.
\end{lemma}

\begin{proof} Let $\psi\in\ss (S^-)$.
By Theorem \ref{1.2}, the function $G_\psi$ is pluriharmonic and
so $\Delta_rG_\psi=0$ and for $j=1,\ldots,r$, $\alpha=1,\ldots,d$,
 we have $\Delta_{jr}^\alpha G_\psi=0$. But
$$
\Delta_{jr}^\alpha G_\psi(s^+)=\int_{S^-}
\Delta_{jr}^\alpha F(s^-s^+)\psi(s^-)\,ds^-,
\quad\Delta_r G_\psi(s^+)=\int_{S^-}
\Delta_rF(s^-s^+)\psi(s^-)\,ds^-
$$
for any $\psi$, so the conclusion follows.
\end{proof}

\begin{proof}[End of the proof of Theorem \ref{mainth}]
We have just proved
\arraycolsep 1pt
\begin{equation}
\begin{array}{ccccl}
\partial_{z_{jr}^\a}\partial_{\overline z_{jr}^\a}F(e)&=&0&\quad&\mbox{for }j\leq r-1,\mbox{ and }
\alpha=1,\ldots,d\\
\partial_{z_{r}}\partial_{\overline z_{r}}F(e)&=&0.&&\\
\end{array}
\end{equation}
Taking any element $k\in$Aut($V$) permuting the chosen Jordan frame and repeating
the above argument, but with the group
$S_k=kSk^{-1}$ instead of $S$, we may prove that $F$
is annihilated by the $S_k$-invariant operators $k\Delta_{jr}^\a k^{-1}$. Therefore,
\arraycolsep 1pt
\begin{equation}
\begin{array}{ccccl}
\partial_{z_{jk}^\a}\partial_{\overline z_{jk}^\a}F(e)&=&0&\quad&\mbox{for }j<k\leq r,\mbox{ and }
\alpha=1,\ldots,d\\
\partial_{z_{k}}\partial_{\overline z_{k}}F(e)&=&0&\quad&\mbox{for }k\leq r\\
\end{array}
\end{equation}
Take any $s\in S$, then $F_s=F\circ s$ satisfies assumptions of
Theorem \ref{mainth}, hence \arraycolsep 1pt
\begin{equation}
\begin{array}{ccccl}
\partial_{z_{jk}^\a}\partial_{\overline z_{jk}^\a}(F\circ s)(e)&=&0&\quad&\mbox{for }j<k\leq r,\mbox{ and }
\alpha=1,\ldots,d\\
\partial_{z_{k}}\partial_{\overline z_{k}}(F\circ s)(e)&=&0&\quad&\mbox{for }k\leq r\\
\end{array}
\end{equation}
which implies (\ref{cond2}) and completes the proof of the main theorem.
\end{proof}

\appendix
\section{A lemma about confluent hypergeometric equations}
\label{appendixa}

\begin{notation}
For $a\in\R$ and $n\in\N$, we denote by $(a)_0=1$, $(a)_1=a$
and $(a)_n=a(a+1)\ldots(a+n-1)$. Further, if $a,c\in\R$ and
$c\notin\Z^-$, then
$$
{}_1F_1(a,c,x)=\sum_{n=0}^{+\infty}\frac{(a)_n}{(c)_n}\frac{x^n}{n!}
$$
defines an entire function (called the {\it confluent
hypergeometric function}) that satisfies the following
differential equation:
$$
xy''+(c-x)y'-ay=0.
$$
Notice that, for $x>0$, if $a>c>0$, then ${}_1F_1(a,c,x)\geq e^x$,
whereas if $a>0>c$, then ${}_1F_1(a,c,x)\geq e^x-P(x)$ where $P$
is a polynomial.
\end{notation}
We are interested in bounded solutions on $]0,+\infty)$ of a
variant of this equation.
\begin{lemma}
\label{hyper} Let $\gamma>0$ and $\beta\geq 0$ and let $k$ be the
smallest integer greater than $\gamma$, that is $k=\ent{\gamma}+1$
if $\gamma$ is not an integer and $k=\gamma$ otherwise. Then, the
equation
\begin{equation}
xg''-\gamma g'-\left(x+\left(\gamma+\beta\right)\right)g=0
\label{hypereq}
\end{equation}
 has one bounded non zero solution $y$ on $(0,+\infty)$ and only one,
up to a constant. For $\beta \neq 0$,
 its
derivatives have finite limits at $0$ up to the order $k$.
Moreover, for every positive integer $N$, there exists two
functions, $a_1$ and $a_2$,  which are smooth up to order $N$ on
$[0,\infty)$, such that $a_1(0)\neq 0$ and
$$
\partial^{k+1}y(x)=
\begin{cases}a_1(x)\ln x+a_2(x) &\mbox{if $\gamma$ is an integer}\\
a_1(x)x^{\gamma-k}+a_2(x) &\mbox{if $\gamma$ is not an integer}\\
\end{cases}.
$$
For $\beta=0$, the solution $y$ is given by $y(x)=e^{-x}$ and is
an entire function.
\end{lemma}

\begin{proof}[Proof]
We leave the computations for the reader when $\beta=0$ since they
are completely elementary.

>From now on, we assume that $\beta>0$. Let us look for solutions
of (\ref{hypereq}) of the form $g(x)=e^{-x}h(2x)$. Then $h$ is a
solution of the equation
\begin{equation}
xh''-(x+\gamma)h'-\left(\frac{\beta}{2}\right)h=0. \label{confeq}
\end{equation}

When $\gamma$ is not a non negative integer, the function
${}_1F_1\left(\frac{\beta}{2},-\gamma,\cdot\right)$ is a solution
of (\ref{confeq}). Moreover, it increases at infinity in such a
way that the corresponding solution of (\ref{hypereq}) is
unbounded. For $\gamma$ a non negative integer, there is still an
entire function which is a solution of the equation (\ref{confeq})
and gives rise to an unbounded solution of (\ref{hypereq}).
Coefficients of its Taylor series vanish up to $-c$. The
uniqueness, up to a constant, of  bounded solutions of
(\ref{hypereq}) on $]0,\infty)$ follows at once since the vector
space of solutions, which is of dimension $2$, can only have a
proper subspace of dimension $1$ of bounded solutions.

Let us now show the existence of a bounded solution on
$]0,+\infty)$ as well as the asymptotic behavior of the $(k+1)$-th
derivative. Using a Laplace transform for (\ref{confeq}), it is
easy to propose another solution for (\ref{hypereq}),
$$
y(x)=e^{-x}\int_0^{\infty}e^{-2xt}
\frac{t^{\frac{\beta}{2}-1}}{(1+t)^{\gamma +\frac{\beta}{2}}}
 \,dt.
$$
It is easily seen that $y$ is a $\cc^\infty$ function on
$(0,+\infty)$,  of class $\cc^k$ on $[0,+\infty)$, and bounded. It
remains to find the explicit expression for the derivative of
order $k+1$ of $y$. Let us first assume that $\gamma$ is not an
integer. We look for another solution of (\ref{hypereq}) which may
be written as $x^{\gamma+1}e^{-x}h(2x)$. Then $h$ is a solution of
$$xh''-(x+\gamma+2)h'-(\gamma+\frac{\beta+1}{2})h=0.$$
So $y$, as a solution of (\ref{hypereq}), may be written
$$
y(x)= \alpha
x^{\gamma+1}e^{-x}{}_1F_1\left(\gamma+\frac{\beta+1}{2},\gamma+2,2x\right)+\tilde\alpha
e^{-x}{}_1F_1\left(\frac{\beta}{2},-\gamma,2x\right).
$$
Moreover $\alpha\neq 0$ since $y$ is bounded. This allows to
conclude in this case, taking derivatives up to order $k+1$.

Let us now consider the case when $\gamma=k-1$ is an integer. We
will use the explicit formula of $y$ as an integral. Moreover,
cutting the integral into two, from $0$ to $1$ and from $1$ to
infinity, we see that the first integral is a $\cc^\infty$
function. So it is sufficient to consider
$$
\int_1^{\infty} e^{-2xt} \frac{t^{\frac{\beta}{2}-1}}{(1+t)^{k
+\frac{\beta}{2}}}
 \,dt =
\int_1^{\infty} e^{-2xt} t^{-k}\left(1+\frac 1t\right)^{k-1
+\frac{\beta}{2}}
 \,dt.
$$
Moreover, taking the Taylor expansion of the term $(1+\frac
1t)^{k-1 +\frac{\beta}{2}}$, it is sufficient to consider terms up
to order $N$, the remaining part giving rise to a $\cc^{N+k}$
function. Finally, after having taken derivatives, it is
sufficient to prove that each term
$$I_j(x):= \int_1^{\infty} e^{-xt} t^{-j-1} \,dt$$
may be written as $\alpha_jx^j \ln x +a_j(x)$, with $\alpha_0\neq
0$ and $a_j$ a $\cc^\infty$ function. Integration by parts allows
to reduce to the case when $j=0$. But a direct computation gives
$$I_0(x)=\int_1^{\infty} e^{-t}\frac {dt}t +\sum_{l=0}^\infty
\frac{(-1)^\ell}{\ell!}\int_x^{1} t^{l-1}
\,dt,$$ from which we conclude easily.

\end{proof}

\section{A lemma about Legendre equations}
\label{appendixb}

\begin{notation} For $c\in\R\setminus\Z^-$, define the entire function
$$
{}_0F_1(c,x)=\sum_{k=0}^\infty\frac{1}{(c)_k}\frac{x^k}{k!}
$$
defines an entire function, called the {\it Legendre function},
that satisfies the following differential equation:
$$
xy''(x)+cy'(x)-y(x)=0.
$$
\end{notation}
For $c$ a negative integer, there is still an entire function
which is a solution of the Legendre equation. Coefficients of  its
Taylor series vanish up to $-c$.

We are interested in bounded solutions of the Legendre equation
for negative values of the parameter.

\begin{lemma}
\label{legen} Let $\beta >0$ and let $k$ be the smallest integer
greater than $\beta$. Then the equation
\begin{equation}
xg''-\beta g'- g=0 \label{legeneq}
\end{equation}
has a bounded non zero solution $y$ on $]0,+\infty)$ whose
derivatives have finite limits at $0$ up to the order $k$. For
every positive integer $N$, there exists two  functions, $a_1$ and
$a_2$,  which are smooth up to order $N$ on $[0,\infty)$, such
that $a_1(0)\neq 0$ and
$$
\partial^{k+1}y(x)=
\begin{cases}a_1(x)\ln x+a_2(x) &\mbox{if $\beta$ is an integer}\\
a_1(x)x^{\beta-k}+a_2(x) &\mbox{if $\beta$ is not an integer}\\
\end{cases}.
$$
Moreover, all non zero  bounded solutions of (\ref{legeneq})
behave as above.
\end{lemma}

\begin{proof}[Proof]
The uniqueness, up to a constant, of  bounded solutions follow
from the fact that the entire functions that are solutions
increase faster than any polynomial at infinity. So the vector
space of solutions, which is of dimension $2$, can only have a
proper subspace of dimension $1$ of bounded solutions. Let us now
show the existence of a bounded solution on $(0,+\infty)$.

This one is given by the  Laplace method:
$$
y(x)=\int_0^{+\infty}\frac{e^{-xt}e^{-1/t}}{t^{2+\beta}}dt.
$$
It is clear that $y$ is rapidly decreasing and has derivatives up
to order $k$ at $0$. It remains to find the explicit formula for
the $k+1$-th derivative. When $\beta$ is not an integer, we look
for another solution of  (\ref{legeneq}) which may be written as
$x^{\beta+1}h(x)$.  Then $h$ is a solution of the equation
$$xh''+(\beta+2)h'-1=0.$$
As a consequence, the function $y$, which is a solution of
(\ref{legeneq}), may be written as $$y(x) = \alpha
x^{\beta+1}{}_0F_1(\beta+2,x)+\tilde\alpha{}_0F_1(-\beta,x),$$
with $\alpha\neq 0$ since ${}_0F_1(-\beta,\cdot)$ is not bounded.
This allows to conclude in this case, taking derivatives up to
order $k+1$.

For $\beta=k$, we use the explicit formula for the function $y$,
that is, up to the constant $(-1)^{k+1}$,
$$
\int_0^1\frac{e^{-xt}e^{-1/t}}{t}\,dt
+\int_1^{+\infty}\frac{e^{-xt}(e^{-1/t}-\sum_{j=0}^{N}(-1/t)^j/j!)}{t}\,dt
+\sum_{j=0}^{N}c_j\int_1^{+\infty}\frac{e^{-xt}}{t^{j+1}}\,dt,
$$
with $c_j=(-1)^j/j!$. Now, the two first integrals are smooth up
to order $N$. It remains to consider each term of the last sum,
for which we conclude as in the previous section.
\end{proof}

\section{Some complements on Harmonic Analysis on the Heisenberg group}

\subsection{A lemma about Laguerre functions}
\label{applaguerre}
We come back to the notations of Section
\ref{subsec:rephei}. Further lemmas we may need can be found e.g.
in \cite{T}.

Recall that we defined a unitary representation $R^\lambda$ of $\hei$ on
$L^2(\R^n)$ by
$$
R^\lambda(\zeta,t)\Phi(x)=e^{2\pi
i\lambda(u.x+u.v/2+t/4)}\Phi(x+v)
$$
The Hermite functions $H_k$ of one variable are
defined by the formula
$$
h_k(x)=(-1)^kc_ke^{x^2/2}\left(\frac{d}{dx}\right)^ke^{-x^2},\quad
k=0,1,2,\ldots,
$$
with $c_k=(\sqrt{\pi}2^kk!)^{1/2}$. For
$k=(k_1,\ldots,k_n)$  a multi-index, we write $h_k=h_{k_1}\ldots
h_{k_n}$.

Now, for $\lambda\not=0$, denote by $h_k^\lambda(x)=
(2\pi\abs{\lambda})^{n/4}h_k\bigl((2\pi\abs{\lambda})^{1/2}x\bigr)$.
Finally, let
$$
e^{\lambda}_\kappa(\omega)=\sum_{|k|=\kappa}(R^{\lambda}(\omega)
h_k^{\lambda},h_k^{\lambda})
$$
and, for $\psi\in C_c^\infty(\R\backslash\{0\})$,
\begin{eqnarray}
\label{5.10}
e^{\psi}_k(\zeta,t)=\int_\R
e_k^\lambda(\zeta,t)\psi(\lambda)\,d\lambda.
\end{eqnarray}
Then, $e^\phi_k\in{\ss}(\C^n\times \R)$.

\begin{lemma}
\label{nor}
Let $\eps>0$.
There exists constants $C,M>0$ such that, for every
$\psi\in C_c^\infty(\R)$ with $\supp\psi\subset\ent{\eps,\eps^{-1}}$,
and every $\kappa\in\N$,
\begin{equation}
\label{norr} \norm{e^\psi_\kappa}_{L^1}\leq C\kappa^M
\bigl(\|\psi\|_{L^{\infty}} +\|\psi'\|_{L^{\infty}}\bigr)
\end{equation}
\end{lemma}

Before proving this lemma we need more information about the function
$\Phi_k=\scal{R^\lambda h_k^\lambda,h_k^\lambda}$.
Let $L_k$ be the $k$-th Laguerre polynomial, i.e.
$$
L_k(t)e^{-t}=\frac{1}{k!}\bigg(\frac{d}{dt}\bigg)^k\left(e^{-t}t^k\right).
$$
Given a multi-index $k=(k_1,...,k_n)$ and
$\zeta\in\C^n$ let
\begin{equation}
\label{5.17}
\mathrm{L}_k(\zeta)=L_{k_1}\left(\frac{1}{2}|\zeta_1|^2\right)\times ...\times
L_{k_n}\left(\frac{1}{2}|\zeta_n|^2\right).
\end{equation}
Then
\begin{equation}
\label{5.18}
\Phi_k(\zeta)=(2\pi)^{-n/2}\mathrm{L}_k(\zeta)e^{-\frac{1}{4}|\zeta|^2},
\end{equation}
(see Formula (1.4.20) \cite{T}).

We will use the following well-known property of the Laguerre
functions (see \cite{BDH}).
\begin{lemma}
\label{5.19}
For every $l,p\in \N$ there exist $c=c(l,p)$ and $M=M(l,p)$ such that,
for every $k\in\N$,
\begin{equation}
\label{5.20}
\int_{0}^{\infty}t^l|\partial_t^pL_k(t)|^2e^{-t}dt\leq ck^M.
\end{equation}
\end{lemma}

\begin{proof}[Proof of lemma \ref{nor}]
To estimate the $L^1(\hei)$ norm of $e^\psi_\kappa$, we first use
Schwartz' inequality to see that it is sufficient to have a bound
for
$$
I_k:=\int_{\hei}(1+t^2)(1+|\zeta_1|^4)\times(1+|\zeta_2|^4)\times\cdots
\times (1+|\zeta_n|^4)\times |\phi_k(\zeta,t)|^2\,dt\,d\zeta,
$$
where
$$
\phi_k(\zeta,t):=\int_{\R}e^{i\lambda t}\psi(\lambda)
\Phi_k(\sqrt{\lambda}\zeta)d\lambda.
$$
We use Parseval Identity in the $t$ variable to write that
$$
I_k=\int_{\C^n}\int_\R(1+|\zeta_1|^4)\times\cdots \times
(1+|\zeta_n|^4)\times
\left(|\psi(\lambda)\Phi_k(\sqrt{\lambda}\zeta)|^2+
|\partial_\lambda\{\psi(\lambda)\Phi_k(\sqrt{\lambda}\zeta)\}|^2\right)
\,d\lambda\,d\zeta.
$$
We find that
$I_k\leq CJ(k) (\|\psi\|_{L^{\infty}}+\|\psi'\|_{L^{\infty}})^2$, with
$$
J(k):=\int_{\C^n}\int_\R(1+|\zeta_1|^4)\times\cdots \times
(1+|\zeta_n|^4)\times
\left(|\Phi_{\alpha,\alpha}(\sqrt{\lambda}\zeta)|^2+
|\partial_\lambda\{\Phi_{\alpha,\alpha}(\sqrt{\lambda}\zeta)\}|^2\right)
\,d\lambda\,d\zeta.
$$
Now  we use Relation (\ref{5.18}), and integrate first in $\zeta$, using Lemma \ref{5.19}.
We then find that $J(k)\leq c\kappa^M$ where $\kappa=|k|$. We have thus proved
(\ref{norr}).
\end{proof}

\subsection{The Fourier Inversion Formula}
\label{heifinv}
In this section, we adapt the Fourier  Inversion Formula
to our choice of representation of the Heisenberg group.
We only give an outline of the proofs, details with a similar choice of representation
may be found in \cite{Fo}, pages 35-37.

The Fourier transform of an integrable
function $f$ on $\hei$ is the operator-valued function
$\lambda\mapsto\hat{f}(\lambda)$ given by
$$
\scal{\hat
f(\lambda)\phi,\psi}=\int_{\hei}\scal{R^\lambda(\w)\phi,\psi}f(\w)d\w
dt.
$$
To avoid any technicality, let us first assume that $f$ is bounded and compactly supported.
One may then show that
\begin{equation}
\label{FIheipre}
\mathrm{tr}\,\bigl(R^\lambda(\zeta,t)\widehat{f}(\lambda)\bigr)
=\abs{\lambda}^{-n}\ff_3f(\zeta,\lambda/4)e^{2i\pi\lambda t/4}
\end{equation}
where $\ff_3$ is the partial Fourier transform in the $t$-variable.
The ordinary Fourier Inversion Formula (in the $t$-variable) then leads to the
Heisenberg Fourier Inversion Formula
\begin{equation}
\label{FIhei}
f(\zeta,t)=c\int_{\R}\mathrm{tr}\,\bigl(R^\lambda(\zeta,t)\widehat{f}(\lambda)\bigr)
\abs{\lambda}^nd\lambda.
\end{equation}
Now, let $\ffi\in\ss(\R)$ be such that $\widehat{\ffi}$ is compactly supported
away from $0$ and apply (\ref{FIhei}) to $f*_3\ffi$, the convolution in the $t$-variable of $f$ and $\ffi$,
to get
$$
f*_3\ffi(\omega)=c\int_{\R}
\mathrm{tr}\,\bigl(R^\lambda(\omega)\widehat{f*_3\ffi}(\lambda)\bigr)
\abs{\lambda}^nd\lambda,
$$
and, using twice (\ref{FIheipre}), this becomes
\begin{equation}
\label{FIhei3}
f*_3\ffi(\omega)=c\int_{\R}
\mathrm{tr}\,\bigl(R^\lambda(\omega)\widehat{f}(\lambda)\bigr)
\widehat{\ffi}(\lambda/4)\abs{\lambda}^nd\lambda.
\end{equation}
Next, rewriting the trace using the Hermite basis leads to
\begin{align}
\mathrm{tr}\,\bigl(R^\lambda(\omega)\widehat{f}(\lambda)\bigr)
=&\sum_{\kappa\in\N}\sum_{k\,:|k|=\kappa}\scal{R^\lambda(\omega)\widehat{f}(\lambda)
h_k^\lambda,h_k^\lambda}\notag\\
=&\sum_{\kappa\in\N}\sum_{k\,:|k|=\kappa}\int_{\hei}
\scal{R^\lambda(\omega)R^\lambda(\eta,s)h_k^\lambda,h_k^\lambda}
f(\w)d\w\notag\\
=&\sum_{\kappa\in\N}\int_{\hei}e_\kappa^\lambda(\omega\w)f(\w)d\w
=\sum_{\kappa\in\N}\int_{\hei}e_\kappa^\lambda(\w)f(\omega^{-1}\w)d\w.\notag
\end{align}
Further, inserting this in (\ref{FIhei3}), inverting summations and
integrations, justified by Lemma \ref{nor} and Fubini's theorem,
we get
\begin{equation}
\label{FIhei2}
f*_3\ffi(\omega)=\sum_{\kappa\in\N}
\int_{\hei}f(\omega^{-1}\w)e_\kappa^{\psi}(\w)d\w,
\end{equation}
where $\psi(\lambda)=\widehat{\ffi}(\lambda/4)|\lambda|^n$.
Finally, note that Identity (\ref{FIhei2}) is also valid for $f\in L^\infty$
since we may apply it to a truncation
$f_R(\omega)=\begin{cases}f(\omega)&\mbox{if }\abs{\omega}\leq
R\\0&\mbox{else}\\ \end{cases}$
and then let $R\to\infty$, again with the help of Lemma \ref{nor}.

\section{Second regularity condition}
\label{appendixd} In this section we show another regularity
condition that leads to the pluriharmonicity of $F$. Since it is
difficult to compare with the previous one, we include it here.
Compared to the previous one, the conditions are given using only
one solvable group $S$ and one special coordinate system.

\medskip
We fix a Jordan frame $c_1,...,c_r$ and the group $S$. Let
\begin{eqnarray*}
\nn_k^-=\bigoplus_{i\le j \le k}V_{ij}\oplus\bigoplus_{i<j\le k}\nn_{ij},&\quad&
\nn_k^+=\bigoplus_{k<i\le j}V_{ij}\oplus\bigoplus_{k<i<j}\nn_{ij},\\
\aa_k^-=\mbox{lin} \{L(c_1),...,L(c_k)\},&&
\aa_k^+=\mbox{lin} \{L(c_{k+1}),...,L(c_r)\}.
\end{eqnarray*}
Then $\nn_k^+$ is an ideal. Clearly
\begin{eqnarray*}
\nn =\nn_k^-\oplus \nn_k^+,&\quad &
\aa =\aa_k^-\oplus \aa_k^+.
\end{eqnarray*}
Let
\begin{eqnarray*}
N_k^- = \exp\nn_k^-,&\quad&
N_k^+ = \exp\nn_k^+,\\
A_k^- = \exp\aa_k^-,&&
A_k^+ = \exp\aa_k^+,\\
S_k^-= N_k^-A_k^-,&&
S_k^+= N_k^+A_k^+.
\end{eqnarray*}
Then
\begin{eqnarray*}
N= N_k^-N_k^+,&\quad&
S= S_k^-S_k^+\\
\end{eqnarray*}
in the sense that
\begin{eqnarray*}
S_k^-\times S_k^+\ni(s^-,s^+)&\mapsto& s^-s^+\in S,\\
N_k^-\times N_k^+\ni(n^-,n^+)&\mapsto& n^-n^+\in N\\
\end{eqnarray*}
are diffeomorphisms.

Let
$$V_k^-=\bigoplus_{i\le j\le k}V_{ij}.$$
$V_k^-$ is the Jordan algebra with the Jordan frame $c_1,...,c_r$
and $S_k^-$ is the solvable group acting simply transitively on $V_k^-+i\Omega_k^-$,
where $\Omega_k^-=$int$\{x^2:\;x\in V_k^-\}$. For every $k=1,...r-1$ we define a subdomain
$$\dd_k=\dd\cap\{z_{lj}=0,z_{jj}=i,\;\; l<j,\;j>k\}.$$
Clearly
$$\dd_k=S_k^-\cdot ie$$
and
$$\partial\dd_k=\partial\dd\cap\{z_{lj}=0,z_{jj}=i,\;\; l<j,\;j>k\}.$$
$\dd_k$ is a symmetric tube domain corresponding to the cone $\Omega_k^-$.
\begin{definition}
Given a tube domain $\dd$ with a fixed Jordan frame $c_1,...,c_r$, the group $S$
and its subgroup $S'$, we say
 that a bounded function is {\it weak regular} if it satisfies
 \eqref{cond1} with $$\Phi(w,b) = \phi(w)+ibc_r$$
 and $k=1,...,\ent{\frac{(r-1)d+1}{2}}+1$.
\end{definition}
 It follows from a theorem by Lassalle \cite{L} that $\phi(\R^{2n-1})$ is
 dense in $\partial\dd$ although smaller than $\widetilde{\partial\dd}$.
 Therefore, weak regularity means certain boundary regularity at any point from a dense subset
 of $\partial\dd$. We are going to assume that the function is weak regular on the domain $\dd$
 and subdomains $\dd_k$. More precisely, given $F$ we write a number of functions $K_{\psi}$
 on the domains $\dd_k$:
 $$K_{\psi}(z) = \int_{S_k^+}F(s^+\cdot z)\psi(s^+)ds^+,$$
 $\psi\in C_C^\8(S_k^+)$, $z\in\dd_k$. $K_{\psi}$ will be identified with its version on the group $S_k^-$:
 $$K_{\psi}(s^-) = \int_{S_k^+}F(s^+s^-\cdot ie)\psi(s^+)ds^+.$$
\begin{definition}
We say that $F$ is regular if $F$ is weak regular and for every $k=1,...r-1$ and $\psi\in C_C^\8(S_k^+)$,
$K_{\psi}$ is weak regular on $\dd_k$ with respect to the Jordan frame $c_1,...c_k$.
\end{definition}
\begin{theorem}
Assume that $F$ is a real, bounded, Hua-harmonic and  regular. Then $F$ is pluriharmonic.
\end{theorem}
\begin{proof}
We are going to prove by induction \eqref{cond2}, which again, in view of \cite{BDH}, implies
pluriharmonicity of $F$. First we notice that weak regularity of $F$ is sufficient to conclude that
$$\Delta_rF=0\quad\mbox{and } \Delta_{jr}^\a F=0$$
(see the proof of Lemma \eqref{lem2.2}). Hence for $j=r$ the conclusion holds.

\medskip
Assume now that $\Delta_{ij}^\a F=0$, $\Delta_{jj}F=0$ for $j>k$ and $i < j$. We consider the operators
$$\hil_j^k=\Delta_{jj}+\frac{1}{2}\Big(
\sum_{\substack{i<j\\ \a}}\Delta_{ij}^\a+\sum_{\substack{j<i\le k\\ \a}}\Delta_{ji}^\a\Big).$$
$\hil_j^k$ have perfect meaning on $S_k^-$ and moreover, we have
\begin{equation}
\label{st}
(\hil_j^kF)(s^+s^-)= (\hil_j^kF_{s^+})(s^-),
\end{equation}
where $F_{s^+}(s^-) = F(s^+s^-)$ is a function on $S_k^-$. In \eqref{st} $\hil_j^k$ on the left-hand
side is considered on $S$, while on the right hand side on $S_k^-$. Finally, $\hil_j^k$
are the strongly diagonal Hua operators on $S_k^-$. By induction the left-hand side of
\eqref{st} vanishes so
$$\hil_j^k K_{\psi}(s^-)=0.$$
Now weak regularity of $K_{\psi}$ implies that
$$\Delta_{ik}^\a K_{\psi}=0,\quad\Delta_{kk}K_{\psi}=0.$$
But
$$\Delta_{ik}^\a K_{\psi}(s^-)=\int_{S_k^+}(\Delta_{ik}^\a F)(s^+s^-)\psi(s^+)ds^+$$
and so $\Delta_{ik}^\a F=0$.
\end{proof}

\end{document}